\documentstyle[11pt]{article}
\newtheorem{theorem}{Theorem}[section]
\newtheorem{lemma}[theorem]{Lemma}
\newtheorem{corollary}[theorem]{Corollary}
\newtheorem{proposition}[theorem]{Proposition}
\newtheorem{rem}[theorem]{Remark}
\newenvironment{remark}{\begin{rem}\rm}{\end{rem}}
\newtheorem{define}[theorem]{Definition}
\newenvironment{definition}{\begin{define}\rm}{\end{define}}
\newtheorem{ex}[theorem]{Example}
\newenvironment{example}{\begin{ex}\rm}{\end{ex}}
\newenvironment{proof}{\par {\em Proof.}}{ \hfill}

\font\twlmsbm=msbm10 scaled \magstep1
\font\egtmsbm=msbm8
\font\sixmsbm=msbm6
\newfam\msbmfam
\textfont\msbmfam=\twlmsbm
\scriptfont\msbmfam=\egtmsbm
\scriptscriptfont\msbmfam=\sixmsbm

\newcommand{\Tr}{\mbox{Tr}\,}
\newcommand{\id}{\mbox{id}}

\newcommand{\Index}{\mbox{Index}\,}

\newcommand{\Hom}{\mbox{Hom}}

\renewcommand{\span}{\mbox{span}}

\newcommand{\Rep}{\mbox{Rep}}
\newcommand{\Corep}{\mbox{Corep}}

\newcommand{\eps}{\varepsilon}
\newcommand{\lact}{\triangleright}
\newcommand{\ract}{\triangleleft}
\newcommand{\la}{\langle\,}
\newcommand{\ra}{\,\rangle}

\renewcommand{\star}{{\dag}}
\newcommand{\rtimes}{{>\!\!\!\triangleleft}}
\newcommand{\ltimes}{{\triangleright\!\!\!<}}
\newcommand{\smrtimes}{>\!\!\triangleleft}

\newcommand{\half}{1/2}

\newcommand{\1}{_{(1)}}
\newcommand{\2}{_{(2)}}

\newcommand{\I}{^{(1)}}
\newcommand{\II}{^{(2)}}

\font\twlmsbm=msbm10 scaled \magstep1
\font\egtmsbm=msbm8
\font\sixmsbm=msbm6
\newfam\msbmfam
\textfont\msbmfam=\twlmsbm
\scriptfont\msbmfam=\egtmsbm
\scriptscriptfont\msbmfam=\sixmsbm
\def\Bbb#1{{\fam\msbmfam\relax#1}}

\title{\bf A GALOIS CORRESPONDENCE FOR II${}_1$ FACTORS AND QUANTUM GROUPOIDS}
\author{Dmitri Nikshych\thanks{UCLA, Department of Mathematics,
405 Hilgard Avenue, Los Angeles, CA 90095-1555; E-mail: nikshych@math.ucla.edu}
\and Leonid Vainerman\thanks{Universit\'e de Strasbourg, D\'epartement de
Math\'ematiques, 7, rue Ren\'e Descartes, F-67084 Strasbourg, France; E-mail:
wain@agrosys.kiev.ua}}

\date{January 4, 2000}
\begin{document}
\maketitle
\vskip 0.5cm
\begin{abstract}

We establish a Galois correspondence for finite quantum groupoid actions
on II${}_1$ factors and show that every finite index and finite depth
subfactor is an intermediate subalgebra of a quantum groupoid crossed
product. Moreover, any such a subfactor is completely and
canonically determined by a quantum groupoid and its coideal $*$-subalgebra.
This allows to express the bimodule category of a subfactor
in terms of the representation category of a corresponding quantum
groupoid and  the principal graph as the Bratteli
diagram of an inclusion of certain $C^*$-algebras related to it.
\end{abstract}


\begin{section}
{Introduction}

This paper continues the research initiated in \cite{NV2}, where
finite index II$_1$-subfactors of depth $2$
were characterized in terms of weak $C^*$-Hopf algebra
crossed products (the latter objects were introduced in \cite{BNSz}).

In what follows, we use the term ``quantum groupoid''
instead of ``weak $C^*$-Hopf algebra'' since we believe it is
important to stress that this algebraic structure provides a natural
non-commutative generalization  of a usual finite groupoid. In particular,
if it is commutative as an algebra (resp.\ co-commutative
as a coalgebra), then it can be identified in a canonical way with
the $C^*$-algebra of functions on a finite groupoid
(resp.\ groupoid algebra). Quantum groupoids also generalize
finite-dimensional Kac algebras (``finite quantum groups'')
\cite{Y}, \cite{NV1}.

According to the characterization obtained in \cite{NV2},
if $N\subset M\subset M_1\subset M_2\subset\dots$
 is the Jones tower constructed from a finite index depth $2$
inclusion $N\subset M$ of II$_1$ factors, then $B=M'\cap M_2$
has a canonical structure of a quantum groupoid
acting outerly on $M_1$ such that $M =M_1^B$ and
$M_2 = M_1\rtimes B$, moreover $A=N'\cap M_1$ is a
quantum groupoid dual to $B$.

In the present paper we extend the above result to show that
quantum groupoids give a uniform description of arbitrary
finite index and finite depth II$_1$-subfactors via
a Galois correspondence. We also
show how to express subfactor invariants such as bimodule
categories and principal graphs in quantum groupoid terms.

After discussing basic definitions and constructions
in Preliminaries (Section 2) we introduce and study coideal
$*$-subalgebras of quantum groupoids (Section 3),
that play an important role in the sequel.

Section 4 starts with a simple observation
(see Proposition~\ref{reducing depth} and
Corollary~\ref{subinclusions})
that any finite depth subfactor $N\subset M$ ($[M:N]< \infty$)
can be viewed as an intermediate for some depth 2 inclusion
$N\subset \widetilde M$. Due to the above characterization result
we have $\widetilde M \cong M\rtimes B$, which allows to describe
$N\subset M$ via a Galois correspondence between intermediate von
Neumann subalgebras of $N\subset \widetilde M$ and left coideal
$*$-subalgebras of $B$ (Theorem~\ref{galois correspondence}).
Thus, every finite depth subfactor is completely determined
by a pair ($B,\,I$), where $I$ is a left coideal $*$-subalgebra
of a quantum groupoid $B$, and can be realized as
$N \subset N\rtimes I$, where $N\rtimes I$ is
a von Neumann algebra generated by $N$ and $I$ inside $N\rtimes B$.
Note that the Galois correspondence for quantum group actions on
factors was established in \cite{E}, \cite{ILP}.

In Section 5 we discuss  an equivalence between the tensor
category of $(N-N)$-bimodules associated with $N\subset M$ and the
co-representation category of $B$ (Theorem~\ref{N-N equivalence}).
Given a quantum groupoid $B$ acting on a II${}_1$ factor $N$
and a pair of its left coideal $*$-subalgebras $H,K$, we define in the spirit
of \cite{Takeuchi} a category $C_{H-K}$ of relative $(B,H-K)$ Hopf
bimodules, whose objects are both $B$-comodules and
$(H-K)$-bimodules such that the $B$-coaction commutes with the
bimodule action and construct a functor from $C_{H-K}$ to the category
of $(N\rtimes H-N\rtimes K)$-bimodules preserving direct sums
and compatible with operations of taking
tensor products and adjoints. In the case when $H$ and $K$
are trivial, $C_{H-K}$ is $\Corep(B)$, the co-representation category
of $B$,  and the above functor is an equivalence. We want to emphasize
that a bimodule category of {\em any} finite depth  subfactor $N\subset M$
(not only depth $2$) is equivalent to $\Corep(B)$ for some $B$.

This functor also allows to express the principal graph
of the inclusion $N\subset M$ in terms of $B$,
as the Bratteli diagram of an inclusion of certain finite-dimensional
$C^*$-algebras related to $B$
(Proposition~\ref{Bratteli diagram = principal graph},
Corollary~\ref{principal graph for depth 2}).

Finally, in Appendix we explicitly write down the structure
maps of a quantum groupoid associated with a finite depth subfactor.

{\bf Acknowledgements.}
The first author would like to thank E.~Effros, S.~Popa,
and M.~Takesaki for their useful comments on the early stages
of this work. The second author is grateful to M.~Enock, V.~Turaev
and J.-M.~Vallin for valuable discussions and to the University Paris-6,
University of Strasbourg, and Katholieke Universiteit Leuven for their kind
hospitality during his work on this paper.

\end{section}


\begin{section}
{Preliminaries}

Throughout this paper we use Sweedler's notation for a
comultiplication, writing $\Delta(b) = b\1 \otimes b\2$.

A {\em weak Hopf  $C^*$-algebra} \cite{BNSz} or {\em quantum groupoid}
$B$  is a finite dimensional $C^*$-algebra with the comultiplication
$\Delta : B\to B\otimes B$, counit $\eps : B\to \Bbb{C}$, and antipode
$S:B\to B$ such that $(B,\Delta, \eps)$
is a coalgebra and the following axioms hold for all $b,c,d\in B$ :
\begin{enumerate}
\item[(1)] $\Delta$ is a (not necessarily unital)  $*$-homomorphism :
$$
\Delta(bc) = \Delta(b)\Delta(c), \quad \Delta(b^*) =
\Delta(b)^*,
$$
\item[(2)] The unit and counit satisfy the identities
\begin{eqnarray*}
\eps(bc\1)\eps(c\2d) &=&\eps(bcd), \\
(\Delta(1)\otimes 1)(1\otimes \Delta(1)) &=& (\Delta\otimes \id)\Delta(1),
\end{eqnarray*}
\item[(3)]
$S$ is an anti-algebra  and anti-coalgebra  map such that
\begin{eqnarray*}
m(\id \otimes S)\Delta(b) &=& (\eps\otimes\id)(\Delta(1)(b\otimes 1)),\\
m(S\otimes \id)\Delta(b) &=& (\id \otimes \eps)((1\otimes b)\Delta(1)),
\end{eqnarray*}
where $m$ denotes the multiplication.
\end{enumerate}
\medskip

The right hand sides of two last formulas define {\em target}
and {\em source counital maps} :
$$
\eps_t(b) = (\eps\otimes\id)(\Delta(1)(b\otimes 1)), \quad
\eps_s(b) = (\id \otimes \eps)((1\otimes b)\Delta(1)),
$$
and play an important role in this theory.

Let us remark that the axiom $(2)$ of the definition of a
quantum groupoid is equivalent to each of the following
axioms expressed in terms of counital maps (\cite{NV1}, \cite{NV2}) :
\begin{enumerate}
\item[(2')]
$\qquad b\eps^t(c) = \eps(b\1c)b\2, \qquad
b\1 \otimes \eps^t(b\2)  = 1\1b\otimes 1\2,$
\item[(2'')]
$\qquad \eps^s(c)b = b\1\eps(cb\2), \qquad
\eps^s(b\1) \otimes b\2 = 1\1 \otimes b1\2,$
\end{enumerate}
These axioms are convenient for concrete computations,
as they show that the properties of the counital maps
$\eps_t$ and $\eps_s$ are similar to those of a counit in
an ordinary Hopf algebra.

The dual vector space $B^*$ has a natural structure of
a quantum groupoid given by dualizing the structure operations of  $B$
(\cite{BNSz}, \cite{NV1}) :
\begin{eqnarray*}
\la \varphi\psi,\, b\ra &=& \la \varphi\otimes\psi,\, \Delta(b)
\ra, \\
\la \Delta(\varphi),\, b\otimes c\ra &=& \la\varphi,\, bc\ra, \\
\la S(\varphi),\, b\ra &=& \la \varphi,\, S(b)\ra, \\
\la \varphi^*,b\ra &=& \overline{ \la \varphi,\, S(b)^*\ra },
\end{eqnarray*}
for all $b,c\in B$ and $\varphi,\psi\in B^*$. The unit of $B^*$ is
$\eps$ and the counit is $1$.

The main difference between finite quantum groupoids and classical
finite-dimensional Hopf $C^*$-algebras (Kac algebras) is that the
images of the counital maps are, in general, non-trivial unital
$C^*$-subalgebras of $B$, called {\em target} and {\em source counital
subalgebras} :
\begin{eqnarray*}
B_t &=& \{b\in B \mid \eps_t(b) =b \}
  =  \{b\in B \mid \Delta(b) =( b \otimes 1)\Delta(1)
  = \Delta(1)( b \otimes 1) \}, \\
B_s &=& \{b\in B \mid \eps_s(b) =b \}
  = \{b\in B \mid \Delta(b) = (1 \otimes b) \Delta(1)
  = \Delta(1)(1 \otimes  b) \}.
\end{eqnarray*}

The counital subalgebras commute elementwise: $[B_t,\, B_s] =0$, we also have
$S\circ\eps^s = \eps^t\circ S$ and $S(B_t) =B_s$. We say that $B$
is {\em connected} \cite{N2} if $B_t \cap Z(B)= \Bbb{C}$ (where $Z(B)$
denotes the center of $B$), i.e.,\ if the inclusion $B_t \subset B$
is connected. $B$ is connected iff $B_t^* \cap B_s^* = \Bbb{C}$
(\cite{N2}, Proposition 3.11). We say that $B$ is {\em biconnected}
if both $B$ and $B^*$ are connected.
\medskip

The antipode of a quantum groupoid is necessarily unique,
invertible, and satisfies $(S\circ *)^2 =\id$. Furthermore, there
exists a canonical positive element $H$ in the center of $B_t$ such that
$S^2$ is an inner automorphism implemented by $G= HS(H)^{-1}$, i.e.,\
$S^2(b) = GbG^{-1}$ for all $b\in B$. The element $G$ is
group-like, i.e.,\ $\Delta(G) =(G\otimes G)\Delta(1)= \Delta(1)(G\otimes G)$.

Quantum groupoids possess integrals in the following sense.
\begin{enumerate}
\item[]
There exists a unique projection $p\in B$,
called a {\em Haar projection}, such that for all $x\in B$ :
$$
x p  = \eps^t(x)p ,\quad S(p) = p, \quad \eps^t(p)=1.
$$
\item[]
There exists a unique positive functional $\phi$ on $B$,
called a {\em normalized Haar functional} (which is a trace iff
$B$ is a weak Kac algebra), such that
$$
(\id\otimes \phi)\Delta = (\eps^t\otimes \phi)\Delta,\quad
\phi\circ S =S,\quad
\phi\circ \eps^t = \eps.
$$
\end{enumerate}

The next proposition establishes a useful invariance property of
the Haar functional.
\begin{proposition}[cf. (\cite{NV1}, 2.3.5)]
\label{invariance}
The normalized Haar functional $\phi$ of $B$ satisfies the
following strong invariance property :
$$
x\1 \phi(yx\2) = S(y\1) \phi(y\2 x), \qquad
\phi(x\1y)x\2 =\phi(xy\1)S(y\2),
$$
for all  $x,y\in B$.
\end{proposition}
\begin{proof}
It follows from the axioms of quantum groupoid that
$$
\eps_t(S(c))b =\eps(cb\1)b\2
$$
for all $b,c\in B$. Using  this identity and
the properties of $\phi$ one computes :
\begin{eqnarray*}
x\1 \phi(yx\2)
&=& 1\1x\1 \phi(y1\2x\2) \\
&=& \eps_s(y\1)x\1 \phi(y\2 x\2) \\
&=& S(y\1) (y\2x)\1 \phi( (y\2x)\2) \\
&=& S(y\1) \eps_t((y\2x)\1) \phi( (y\2x)\2) \\
&=& S(y\1) 1\2 \phi( \eps(1\1 (y\2x)\1) (y\2x)\2 ) \\
&=& S(y\1) 1\2 \phi (\eps_t(S(1\1)) y\2 x ) \\
&=& S(y\1) S(1\1) \phi( 1\2y\2 x) = S(y\1) \phi(y\2 x).
\end{eqnarray*}
The second identity is similar.
\end{proof}
\medskip

When $B$ is connected, there exists a unique (non-degenerate)
Markov trace (\cite{GHJ}) $\tau$ for the inclusion $B_t\subset B$
normalized by $\tau(1)=\dim B_t$.
This trace is related to the Haar functional $\phi$ by $\phi(x) =
\tau(HS(H)x)$ (\cite{NV2}, 5.7), where $H$ is the canonical central positive
element  in $B_t$ described earlier.

\begin{corollary}
\label{tau-invariance}
For all $x,y\in B$ we have
\begin{eqnarray*}
x\1 \tau(yx\2) &=& S(y\1) G \tau(y\2 x),\\
\tau(x\1y)x\2 &=&  \tau(xy\1)G^{-1}S(y\2),
\end{eqnarray*}
where $G= HS(H)^{-1}$ is the canonical  element implementing $S^2$.
\end{corollary}
\begin{proof}
From Proposition~\ref{invariance} we have
$$
x\1 \tau(HS(H)yx\2) = S(y\1) \tau(HS(H)y\2 x),
$$
and replacing $y$ by $HS(H)y$ we get the first identity,
the second one is similar.
\end{proof}
\medskip

The following notions of action, crossed product, and fixed point
subalgebra were introduced in \cite{NSzW}.
A (left) {\em action} of a quantum groupoid $B$ on a
von Neumann algebra $M$ is a linear map
$$
B\otimes M \ni b\otimes x \mapsto (b\lact x)\in M
$$
making $M$ into a left $B$-module such that for all $b\in B$ the
map
$b\otimes x \mapsto (b\lact x)$ is weakly continuous and
\begin{enumerate}
\item[(1)]
$b\lact xy = (b\1\lact x)(b\2\lact y),$
\item[(2)]
$(b\lact x)^* = S(b)^* \lact x^*, $
\item[(3)]
$b\lact 1 = \eps^t(b)\lact 1$, and  $b\lact 1 =0$ iff
$\eps^t(b)=0.$
\end{enumerate}
\medskip

A {\em crossed product} algebra
$M\rtimes B$ is constructed on the relative tensor product
$M \otimes_{B_t} B$, where  $B$ is a left $B_t$-module
via multiplication and $M$ is a right $B_t$-module via
$x\ract z = S(z)\lact x = x(z\lact 1)$.
Let $[x\otimes b]$  denote the class of $x\otimes b$ in $M\rtimes B$.
A $*$-algebra structure on $M\rtimes B$ is defined by
$$
[x\otimes b][y\otimes c] = [x(b\1\lact y) \otimes b\2 c], \quad
[x\otimes b]^* = [(b\1^*\lact x^*) \otimes b\2^*] \\
$$
for all $x,y\in M,\,b,c \in B$. It is possible to show that this
abstractly defined $*$-algebra $M\rtimes B$ is *-isomorphic to a
weakly closed algebra of operators on some Hilbert space
\cite{NSzW}, i.e., $M\rtimes B$ is a von Neumann algebra.

The collection
$ M^B =\{ x\in M \mid b\lact x = \eps^t(b)\lact x,\ \forall b\in B\}$
is a von Neumann subalgebra of $M$, called a {\em fixed point subalgebra}.
The relative commutant $M^\prime \cap M \rtimes B$ always contains
a *-subalgebra isomorphic to $B_s$ \cite{NV2}. The action of $B$ is
called minimal if $B_s = M^\prime \cap M\rtimes B$.

It was shown in \cite{NV2} that finite index depth $2$ subfactors
of II${}_1$-factors can be characterized in terms of quantum groupoids.
Namely, if $N\subset M$
is such a subfactor ($[M:N]= \lambda^{-1}$) and
$$
N \subset M \subset M_1 \subset M_2 \subset \cdots
$$
is the corresponding Jones tower, $M_1 = \langle M,\,e_1 \rangle,\,
M_2 = \langle M_1,\,e_2 \rangle,\dots,$ where $e_1 \in N^\prime
\cap M_1,\,e_2 \in M^\prime \cap M_2,\cdots$
are the Jones  projections. The depth 2 condition means that
$N^\prime \cap M_{2}$ is the basic construction of the inclusion
$N^\prime \cap M  \subset  N^\prime \cap M_{1}$. Let $\tau$ be
the trace on $M_2$ normalized by $\tau(1)=1$. There is a canonical
non-degenerate duality form between $A=N'\cap M_1$ and $B=M'\cap
M_2$ defined by
$$
\la a,\, b \ra = \lambda^{-2}\tau(ae_2e_1Hb),
$$
for all $a\in A$ and $b\in B$,
where $H$ is a central element in $M'\cap M_1$ canonically defined
by the property $\tau(Hz)=\Tr(z)$, where $z\in M'\cap M_1$ and $\Tr$
is the trace of the regular representation of $M'\cap M_1$ on
itself (in other words, $H$ is the index \cite{W1}
of $\tau\vert_{M'\cap M_1}$).

Using this duality, one defines the comultiplication, the counit and
the antipode of $B$ as follows :
\begin{eqnarray*}
\la a_1\otimes a_2,\, \Delta(b) \ra &=&   \la a_1a_2,\, b \ra, \\
\eps(b) &=& \la 1,\, b\ra =\lambda^{-1}\tau(e_2Hb), \\
 S(b) &=& J(HbH^{-1})^*J,
\end{eqnarray*}
for all $a,a_1,a_2\in A$ and $b\in B$,
where $J$ is the canonical modular involution on $L^2(M_1)$
and $b\mapsto Jb^*J$ is a $*$-anti-automorphism of $B=M'\cap M_2$.
The above expression for $S$ follows from the explicit formula
(\cite{NV2}, 4.5(i)).

With these operations and involution $b^\star = S(H)^{-1}b^*S(H)$, the
$*$-algebra $B$ becomes a biconnected quantum groupoid
and  $A$ becomes its dual (see \cite{NV2} for the proof).

The counital subalgebras of $B$ are $B_s = M_1'\cap M_2$ and
$B_t = M'\cap M_1$, moreover $H$ is the canonical element of $B_t$.
The map
$$
\lact : B\otimes M_1\to M_1: b\otimes x\mapsto \lambda^{-1}E_{M_1}(bxe_2)
$$
defines a left action of $B$ on $M_1$, such that $M=M_1^B$ is the
fixed point subalgebra
for this action (here $E_{M_1}$ denotes the  $\tau$-preserving
conditional expectation on $M_1$) and
$$
\theta : M_1\rtimes B \to M_2 : [x\otimes b]
\mapsto x S(H)^{1/2}bS(H)^{-1/2}
$$
is an isomorphism of von Neumann algebras.

It is straightforward to show that the above left action of $B$ on $M_1$
extends to an action $b\lact\xi$ of $B$ on $L^2(M_1):=L^2(M_1,\tau)$ such that
$$
(b\lact\xi,\eta)=(\xi,S(H)b^*S(H)^{-1}\lact \eta),
$$
where $b\in B,\xi,\eta\in L^2(M_1)$ (recall that the involution in
$B=M_2\cap M'$  is different from the one in $M_2$ - see above).
This means  that $L^2(M_1)$ equipped with a scalar product
$$
(\xi,\eta)_{\widetilde{L^2(M_1)}}=(S(H)\lact \xi,\eta)_{L^2(M_1)}
$$
is a unitary left $B$-module.

One can also introduce a right action of $B$ on $L^2(M_1)$ by
setting $\xi\ract b=S^{-1}(b)\lact\xi$. This makes $L^2(M_1)$ a
$(B_t-B_t)$- and $(B_s-B_s)$-bimodule (here and in what follows, the term
bimodule means a unitary bimodule, see, e.g., \cite{JS}).
Since $L^2(M_1)$ is also an $(M_1-M_1)$-bimodule with respect to left and right
multiplications, the properties of the left action of $B$ on $M_1$ give:
$$
b\lact(a\xi)=(b_{(1)}\lact a)(b_{(2)}\lact\xi),\qquad
b\lact(\xi a)=(b_{(1)}\lact\xi)(b_{(2)}\lact a),
$$
$$
J(b\lact\xi)=S(b)^*\lact J\xi:=J\xi\ract b^*,\qquad
(J\xi,J\eta)_{\widetilde{L^2(M_1)}}=(G\lact \eta,\xi)_{\widetilde{L^2(M_1)}},
$$
where $a\in M_1,\xi\in L^2(M_1)$ and $J:a\mapsto a^*$
is the canonical modular involution on $L^2(M_1)$.
\end{section}


\begin{section}
{Coideal $*$-subalgebras}
\begin{definition}
A  left (resp.\ right) {\em coideal} of a quantum groupoid
$B$ is a linear subspace $I\subset B$ such that
$\Delta(I) \subset B\otimes I$
(resp.\ $\Delta(I) \subset I\otimes B$ ).
A left (resp.\ right) {\em coideal
$*$-subalgebra} is a unital $C^*$-subalgebra $I\subset B$ which
is a left (resp.\ right) coideal.
\end{definition}

Note that the target counital subalgebra $B_t$
(resp.\ the source counital subalgebra $B_s$)
is a left (resp.\ right) coideal $*$-subalgebra of $B$
contained in every left (resp.\ right) coideal $*$-subalgebra
$I\subset B$.
It is easy to check that a linear subspace
$I\subset B$ is a left (resp.\ right) coideal iff
it is invariant under the right (resp.\ left) dual action of $B^*$
iff
its {\em annihilator} $I^0 = \{ a\in B^* \mid
\la a,\, b\ra =0,\ \forall b\in I\} \subset B^*$ is a left
(resp.\ right) ideal in $B^*$ iff $S(I)$ is a right (resp.\ left)
coideal.
Note that if $I$ is a left coideal $*$-subalgebra, then
$u^*Iu$ is a left coideal $*$-subalgebra for any unitary $u\in
B_s$. In particular, there can be infinitely many non-equal
conjugated coideal $*$-subalgebras.

For any quantum groupoid $B$, the set $\ell(B)$ of left
coideal $*$-subalgebras is a lattice under the
usual operations:
$$
I_1 \wedge I_2 = I_1\cap I_2, \qquad I_1 \vee I_2 = (I_1\cup I_2)''
$$
for all $I_1,I_2\in \ell(B)$. The smallest element of $\ell(B)$ is $B_t$
and the greatest element is $B$.

\begin{proposition}
\label{tilde I}
If $I\subset B$ is a left coideal $*$-subalgebra of $B$, then
$\tilde{I} = G^{-\half} S(I) G^{\half}$ is a right coideal
$*$-subalgebra
of $B$. The map $I\mapsto \tilde{I}$ is an isomorphism of lattices.
\end{proposition}
\begin{proof}
Clearly, $\tilde{I}$ is a subalgebra. Let $c\in \tilde{I}$,
$c= G^{-\half} S(b) G^{\half}$ for some $b\in I$.
Then, using the group-like propery of $G$,  we have :
\begin{eqnarray*}
c^* &=&
G^{\half} S^{-1}(b^*) G^{-\half} = G^{-\half} S(b^*) G^{\half} \in
\tilde{I},\\
\Delta(c) &=&
G^{-\half} S(b\2) G^{\half} \otimes G^{-\half} S(b\1) G^{\half} \in
\tilde{I}
\otimes B,
\end{eqnarray*}
therefore $\tilde{I}$ is a $*$-invariant right coideal. It is easy
to see that the map $I\mapsto \tilde{I}$ preserves the lattice structure.
\end{proof}
\medskip

We will show that $\ell(B)$ is the dual lattice
of $\ell(B^*)$, i.e., $\ell(B) = \breve{\ell}(B^*)$. The following
proposition describes an explicit isomorphism between these lattices.

\begin{proposition}[cf.\ (\cite{ILP}, 4.6)]
\label{dual lattices}
Let $T\subset B$ be a selfadjoint subset and $I$ be the minimal
right coideal $*$-subalgebra of $B$ containing $T$.
Then $T'\cap B^* \subset B^*\rtimes B$ is a left coideal $*$-subalgebra
of $B^*$ and  $T'\cap B^* = I'\cap B^*$.

If we denote this coideal subalgebra by $I^d$, then the map
$\delta : I \mapsto \tilde{I}^d$ defines a lattice anti-isomorphism
between  $\ell(B)$ and $\ell(B^*)$.
\end{proposition}
\begin{proof}
Obviously, $T'\cap B^*$ is a $*$-subalgebra of $B^*$. In order to
prove that it is a left coideal, we need to show that it remains invariant
under the right dual action of $B$, i.e., that $(x\ract a)$ belongs to
$T'$ for all $a\in B,~x\in T'$. The latter means that
$[x\1\otimes  (t\ract x\2)] = [x\otimes t]$ for all $t\in T$.
Applying $a\in B$ to the above identity on the left
(i.e., using the right dual action of $B$ on $B^*\ltimes B =
B^*\rtimes B$), we get
$$
[(x\ract a)\1 \otimes (t\ract(x\ract a)\2)] =
[(x\1\ract a) \otimes (t \ract x\2)] = [ (x\ract a) \otimes t]
$$
therefore, $(x\ract a)\in T'\cap B^*$ for all $a\in B$.
Note that $I$ is generated, as an algebra, by elements of the form
$(y \lact t)$, with $t\in T$ and $y\in B^*$. We need to show that
$I'\cap B^* \subset T'\cap B^*$, i.e., that any $x$ from $T'\cap
B^*$ commutes with $(y \lact t)$. The latter follows from considering
the left dual action of $B^*$ on $B^*\rtimes B$ :
$$
[(y\lact t)\1 \lact x \otimes (y\lact t)\2] =
[(t\1\lact x) \otimes (y\lact t\2)] = [x\otimes (y\lact t)].
$$
The opposite inclusion is obvious.

Since $I_1'\cap I_2' = (I_1\vee I_2)'$ for all $I_1,\,
I_2\in \ell(B)$, the map $\delta : I\mapsto \tilde{I}^d$
is a homomorphism of lattices. Its inverse is given
by the composition of the maps
$\delta(I) \mapsto \delta(I)'\cap B \subset B^* \ltimes B$
and $I\mapsto \tilde{I}$,  since we have
$$
\delta(I)'\cap B  = (\tilde{I}'\cap B^*)'\cap B =
(\tilde{I}\vee (B^*)')\cap B = \tilde{I}\vee ((B^*)'\cap B)
=\tilde{I},
$$
for all $I\in \ell(B)$. Therefore, $\delta$ is an isomorphism.
\end{proof}

\begin{definition}
A left coideal $*$-subalgebra $I\subset B$ is said to be
{\em connected} if $ Z(I)\cap B_s =\Bbb{C}$.
\end{definition}
\medskip
To justify this definition, note that if $I=B$, then this is
precisely the definition of $B$ being connected, and if $I=B_t$, then this
definition is equivalent to $B^*$ being connected (\cite{N2}, 3.10, 3.11).

Let $I\subset B$ be a connected left coideal $*$-subalgebra of $B$,
then there is a uniquely determined positive element $x_I\in I$
such that $ \eps(b) =\tau(x_Ib)$, for all $b\in I$.

\begin{proposition}
\label{x_I}
For any system $\{f^\alpha_{rs}\}$ of matrix units in
$I=\Sigma_\alpha\,
M_{n_\alpha}(C)$ the value of the comultiplication on $x_I$ is
$$
\Delta(x_I) = \sum_{\alpha rs}\,\frac{1}{\tau(f^\alpha_{ss})}\,
              G S^{-1}(f^\alpha_{sr}) \otimes f^\alpha_{rs}
 = \sum_{\alpha rs}\,\frac{1}{\tau(f^\alpha_{ss})}\,
              S(f^\alpha_{sr}) G \otimes f^\alpha_{rs}.
$$
We also have $S(x_I) =x_I G^{-1}$.
\end{proposition}
\begin{proof}
From Corollary~\ref{tau-invariance} we get
\begin{eqnarray*}
b
&=& b\1\eps(b\2) = b\1\tau(x_Ib\2) \\
&=& S(H^{-1}{x_I}\1) \tau(H{x_I}\2b).
\end{eqnarray*}
If we write $\Delta(x_I) =\Sigma_{\alpha rs}\, g^\alpha_{sr}
\otimes
f^\alpha_{sr}$ with $g^\alpha_{sr} \in B$, then applying the above
identity to $b= f^\beta_{kl} H^{-1}$ on has $f^\beta_{kl}H^{-1} =
\tau(f^\beta_{ll}) S(H^{-1}g^\beta_{lk})$, therefore
$ g^\beta_{lk} = \frac{1}{\tau(f^\beta_{ll})}\, G
S^{-1}(f^\beta_{kl})$.
Comparing $\Delta(S(x_I))$ and $\Delta(x_I)$, we get the last
identity.
\end{proof}

\begin{definition}
\label{distinguished projection}
Let $e_I\in I$ be the support of the restriction of $\eps$ on $I$,
i.e., $e_I=e$, where $e$ is the minimal projection having property
$\eps(ebe)=\eps(b)$ for all $b\in I$ (note that $\eps$ is faithful
on
$e_I I e_I$). We will call $e_I$ a {\em distinguished projection}
of $I$.
\end{definition}

Note that $x_Ie_I = e_Ix_I = x_I$ and $e_I$ is the minimal
projection with this property, i.e., $e_I$ is the support of $x_I$.
Also, it is easy to see that $I_1\subset I_2$ implies $e_{I_2} \leq e_{I_1}$.

\begin{proposition}
\label{Haar}
Let $I\subset B$ be a left coideal $*$-subalgebra. Then the
distinguished projection $e_I$ satisfies the following Haar
property :
$$
be_I = \eps_t(b)e_I,\qquad \mbox{for all}~b\in I.
$$
\end{proposition}
\begin{proof}
Since $\eps(xy) = \eps(x\eps_t(y))$ for all $x,y\in B$, we get
$$
\eps( (\eps_t(b)-b)^* (\eps_t(b)-b)) =0,\qquad \forall~b\in B,
$$
which implies
$\eps( e_I (\eps_t(b)-b)^* (\eps_t(b)-b)e_I)) =0,\qquad \forall~b\in I$.

Therefore, $(\eps_t(b)-b)e_I =0$, since $\eps\vert_I(e_I\cdot e_I)$ is
faithful on $e_I I e_I$.
\end{proof}

\begin{remark}
\begin{enumerate}
\item[(i)]
For right coideal $*$-subalgebras one can prove a similar identity
$e_Ib=e_I\eps_s(b)$.
\item[(ii)]
For $I=B$ Proposition~\ref{Haar} is the Haar theorem for quantum groupoids (\cite{NV1}, 2.2.5), (\cite{BNSz}, 4.5),
$e_B$ is the Haar projection, and $x_B$ is a scalar multiple of
$e_B$ (since $e_B$ is minimal in $B$ (\cite{BNSz}, 4.6)).
\item[(iii)]
For $I = B_t$ one has $x_{B_t}=H$ and $e_{B_t} =1$.
\end{enumerate}
\end{remark}

\begin{corollary}
\label{expectation EI}
The map $E_I(y) =y\1\tau(x_Iy\2),~y\in B$ is the
$\tau$-preserving conditional expectation from $B$ to $I$.
\end{corollary}
\begin{proof}
Using the relation from Proposition~\ref{tau-invariance},
the formula of Proposition~\ref{x_I}, and that $S(G)=G^{-1}$ we get
\begin{eqnarray*}
\tau(b E_I(y))
&=& \tau(y\1 b)\tau(x_Iy\2) \\
&=& \tau(yb\1)\tau(x_IG^{-1}S(b\2))\\
&=& \tau(yb\1)\tau(x_Ib\2) = \tau(yb\1)\eps(b\2) = \tau(yb),
\end{eqnarray*}
therefore $E_I$ is the $\tau$-preserving
conditional expectation from $B$ to $I$.
\end{proof}

\begin{remark}
\label{Markov cond.exp}
From the explicit form of the isomorphism $\theta:M_1\rtimes B\to M_2$
one can see that the map
$[x\otimes b]\mapsto [x\otimes E_{M_1}(S(H)^{1/2}bS(H)^{-1/2})]$
from $M_1\rtimes B$ onto $M_1$ is the image of $E_{M_1}:M_2\to M_1$
under $\theta^{-1}$, and thus it is a trace preserving
projection in $M_1\rtimes B$ onto $M_1$. Note that the map
$E_{M_1}: B\to B_t$ is uniquely defined by the relation
$\tau(zb)=\tau(zE_{M_1}(b))\ (\forall z\in B_t, b\in B)$,
and the same relation determines the $\tau$-preserving conditional
expectation $E_{B_t}$ from Corollary \ref{expectation EI}.
So, these two maps  coincide, and the formula for
$E_{B_t}(x) =x\1 \otimes\tau(Hx\2)$ shows that
$E_{M_1}(S(H)^{1/2}b S(H)^{-1/2})=E_{M_1}(b)$.

As a result, the above inverse image of $E_{M_1}:M_2\to M_1$ is the
map $[x\otimes b]\mapsto[x\otimes E_{B_t}(b)]$.
Together with Corollary \ref{expectation EI}
this gives the following expression for $\tau_{M_1\smrtimes B}$:
\begin{eqnarray*}
\tau_{M_1\smrtimes B}([x\otimes b])
&=& \tau_{M_1}(x(E_{B_t}(b)\lact 1))
 = \tau_{M_1}(x(b_{(1)}\lact 1))\tau(Hb_{(2)})\\
&=& \tau_{M_1}(x(S(H1_{(1)})G\lact 1))\tau(1_{(2)}b)
 =  \tau_{M_1}(x(H\lact 1)\tau(b) \\
&=& \tau_{M_1}(S(H)\lact x)\tau(b)=\tau_{M_1}(x)\tau(Hb),
\end{eqnarray*}
where $x\in M_1,b\in B$. Since $\tau_{M_1}(H\lact 1)=\tau(H)=1$,
we have
$$
\tau_{M_1\smrtimes B}([x\otimes 1])=\tau_{M_1}(x)
\qquad  \mbox{ and } \qquad
\tau_{M_1\smrtimes B}([1\otimes b])=\tau(b).
$$
Then one can write down the GNS inner product on $M_1\rtimes B$ as
$$
([x\otimes b],[y\otimes c])_{M_1\smrtimes B}
=(x,y)_{\widetilde{L^2(M_1)}}(b,c)_{B},
$$
where $(\cdot,\cdot)_{B}$ is the GNS-scalar product
on $B$ with respect to the Markov trace.
\end{remark}
\end{section}


\begin{section}
{A Galois Correspondence}

Let $N\subset M$ be a finite depth inclusion of II${}_1$ factors
with finite index $\lambda^{-1}=[M:N]$ and
$$
N\subset M\subset M_1\subset M_2\subset\cdots
$$
be the corresponding Jones tower, $M_i=\la M_{i-1},\,e_i\ra$,
where $e_i\in M_{i-2}'\cap M_i,\ i=1,2,\dots$ are Jones projections
(we denote $M_{-1}=N$ and $M_0 =M$).
Let $n$ be the depth (\cite{GHJ}, 4.6.4) of $N\subset M$, i.e.,
\begin{eqnarray*}
n &=& \min\{ k\in \Bbb{Z}^+ \mid \dim Z(N'\cap M_{k-2}) = \dim
Z(N'\cap M_k) \}.
\end{eqnarray*}

The case $n=2$ is completely understood in \cite{NV2}, where it was shown
that the symmetries of depth 2 subfactors are described by quantum groupoids (see
Preliminaries). For the case of general depth $n\geq 2$ we have the following
result.

\begin{proposition}
\label{reducing depth}
For all $k\geq 0$ the inclusion $N\subset M_k$ has depth
$d+1$, where $d$ is the smallest positive integer $\geq
\frac{n-1}{k+1}$.
In particular, $N\subset M_i$ has depth $2$ for all $i\geq n-2$.
\end{proposition}
\begin{proof}
Note that $\dim Z(N'\cap M_i) = \dim Z(N'\cap M_{i+2})$ for all
$i\geq n-2$.
By \cite{PP2}, the tower of basic construction for $N\subset M_k$
is
$$
N\subset M_k \subset M_{2k+1} \subset M_{3k+2} \subset \cdots,
$$
therefore, the depth of this inclusion is equal to $d+1$, where $d$
is the smallest positive integer such that $d(k+1)-1 \geq n-2$.
\end{proof}

\begin{corollary}
\label{subinclusions}
Any finite depth subfactor $N\subset M$ is an intermediate
subfactor of some depth $2$ inclusion.
\end{corollary}
\begin{proof}
Consider $N\subset M \subset M_{k},~k\geq n-2$.
\end{proof}
\medskip

The last result means that $N\subset M$ can be realized as an intermediate
subfactor of a crossed product inclusion $N\subset N\rtimes B$ for some quantum groupoid $B$ :
$$
N\subset M \subset N\rtimes B.
$$
Recall that in the case of a usual $C^*$-Hopf algebra (i.e., Kac algebra)
action there is a Galois correspondence between intermediate von Neumann
subalgebras of $N\subset N\rtimes B$ and left coideal $*$-subalgebras of
$B$ \cite{ILP},\cite{E}. Thus, it is natural to ask about
a quantum groupoid analogue of this correspondence.

Clearly, the set $\ell(M_1\subset M_2)$ of intermediate von Neumann
subalgebras of $M_1\subset M_2$ forms a lattice under the
operations
$$
K_1 \wedge K_2 = K_1\cap K_2, \qquad K_1 \vee K_2 = (K_1\cup K_2)''
$$
for all $K_1,K_2\in \ell(M_1\subset M_2)$. The smallest element
of this lattice is $M_1$ and the greatest element is $M_2$.

Given a left (resp.\ right)  action of $B$ on a von Neumann algebra
$N$, we will denote (by an abuse of notation)
$$
N\rtimes I = \span\{ [x\otimes b] \mid x\in N,\, b\in I\} \subset
N\rtimes B.
$$

The next theorem establishes a Galois correspondence between
intermediate von Neumann subalgebras of depth $2$ inclusions of
II${}_1$ factors and coideal $*$-subalgebras of a quantum groupoid,
i.e., a lattice isomorphism between $\ell(M_1\subset M_2)$ and $\ell(B)$.

\begin{theorem}
\label{galois correspondence}
Let $N\subset M\subset M_1\subset M_2\subset\cdots$ be the tower
constructed from a depth $2$ subfactor $N\subset M$, $B=M'\cap M_2$
be the corresponding quantum groupoid, and $\theta$ be
the isomorphism between $M_1\rtimes B$ and $M_2$ (\cite{NV2}, 6.3).
Then the following formulas
\begin{eqnarray*}
\phi &:& \ell(M_1\subset M_2) \to \ell(B) :
     K\mapsto \theta^{-1}(M'\cap K) \subset B\\
\psi &:& \ell(B) \to \ell(M_1\subset M_2) :
     I \mapsto \theta(M_1\rtimes I)\subset M_2.
\end{eqnarray*}
define isomorphisms between $\ell(M_1\subset M_2)$ and $\ell(B)$
inverse to each other.
\end{theorem}
\begin{proof}
First, we need to check that $\phi$ and $\psi$ are indeed maps
between the specified lattices. It follows immediately from the
definition of the crossed product that $M_1\rtimes I$ is a von
Neumann subalgebra of $M_1\rtimes B$, therefore $\theta(M_1\rtimes I)$
is a von Neumann subalgebra of $M_2=\theta(M_1\rtimes B)$,
so $\psi$ is a map to $\ell(M_1\subset M_2)$.
To show that $\phi$ maps to $\ell(B)$,
it is enough to show that the annihilator $(M'\cap K)^0 \subset B^*$
of $ M'\cap K \subset B$ is a left ideal in $B^*$.

For all $x\in A,\, y\in (M'\cap K)^0$, and $b\in M'\cap K$ we have
\begin{eqnarray*}
\la xy,\, b\ra
&=& \lambda^{-2}\tau(xye_2e_1Hb) = \lambda^{-2}\tau(ye_2e_1Hbx)\\
&=& \lambda^{-3} \tau(ye_2e_1E_{M'}(e_1Hbx))
  = \la y,\,\lambda^{-1}E_{M'}(e_1Hbx)\ra,
\end{eqnarray*}
and it remains to show that $E_{M'}(e_1Hbx) \in M'\cap K$.
By (\cite{GHJ}, 4.2.7), the square
$$
\begin{array}{ccc}
K & \subset & M_2 \\
\cup &      & \cup \\
M'\cap K &\subset & M'\cap M_2
\end{array}
$$
is commuting, therefore, $E_{M'}(K) \subset M'\cap K$.
Since $e_1Hbx\in K$, we have $xy \in (M'\cap K)^0$, i.e.,
$(M'\cap K)^0$ is a left ideal and $\phi(K)= \theta^{-1}(M'\cap K)$
is a left coideal $*$-subalgebra.

Clearly, $\phi$ and $\psi$ preserve $\wedge$ and $\vee$,
moreover $\phi(M_1) =B_t,\,\phi(M_2)=B$ and $\psi(B_t)=M_1,\,
\psi(B) = M_2$, therefore they are morphisms of lattices.

To see that they are inverses for each other, we first observe that
the condition $\psi\circ\phi =\id$ is equivalent to $M_1(M'\cap K)= K$,
and the latter follows from applying the conditional expectation $E_K$
to $M_1(M'\cap M_2)= M_1B = M_2$. The condition $\phi\circ\psi =\id$
translates into $\theta(I) = M'\cap \theta(M_1 \rtimes I)$.
If $b\in I,\, x\in M=M_1^B$, then
\begin{eqnarray*}
\theta(b)x
&=& \theta([1\otimes b][x\otimes 1])
     = \theta([(b\1\lact x) \otimes b\2]) \\
&=& \theta([x(1\2\lact 1)\otimes \eps(1\1 b\1)b\2 ])
     = \theta([x\otimes b]) =x\theta(b),
\end{eqnarray*}
i.e., $\theta(I)$ commutes with $M$. Conversely, if
$x\in M'\cap \theta(M_1 \rtimes I)\subset B$, then $x=\theta(y)$ for some
$y\in (M_1 \rtimes I) \cap B = I$, therefore $x\in \theta(I)$.
\end{proof}
\medskip

The following proposition describes the center of $K= M_1\rtimes I$
and the first relative commutant in terms of $I$.
\begin{proposition}
\label{center and relative commutant}
In the above situation,
\begin{enumerate}
\item[(i)]
$ Z(K) =  Z(M_1\rtimes I) =  Z(I)\cap B_s$,
\item[(ii)]
$ M_1'\cap K = M_1'\cap M_1\rtimes I = I \cap B_s$.
\end{enumerate}
\end{proposition}
\begin{proof}
Recall that $B_s =M_1'\cap M_2$. If $x\in  Z(I)\cap B_s \subset K$,
then
$x$ commutes with both $I$ and $M_1$ and, therefore, with
$K=M_1\rtimes I$,
i.e., $x\in  Z(K)$.
Conversely, if $x\in  Z(K)$ then $x\in K'\cap K\subset M_1'\cap M_2
=B_s$
and $x\in M'\cap K$, so $x\in  Z(M'\cap K) = Z(I)$ and (i) follows.
To prove (ii), note that
since $B_s = M_1'\cap M_2\subset M_1'$ and $I=M'\cap K\subset K$,
we have
$ M_1'\cap K \subset (M_1'\cap M_2)\cap (M'\cap K) = B_s \cap I
\subset M_1'\cap K$.
\end{proof}

\begin{corollary}
\begin{enumerate}
\item[(i)]
$K = M\rtimes I$ is a factor iff $ Z(I)\cap B_s =\Bbb{C}$.
\item[(ii)]
The inclusion $M_1 \subset K = M_1\rtimes I$ is irreducible iff
$B_s \cap I=\Bbb{C}$.
\end{enumerate}
\end{corollary}

\begin{corollary}
There is no subfactor $N\subset M$ of depth $n$ and index
$\root{k} \of {p}$, where $p$ is prime and  $k\geq n-1$
(unless $n=2$ and $k=1$, in which case $M\cong N\rtimes
\Bbb{Z}/p\Bbb{Z}$).
\end{corollary}
\begin{proof}
Suppose that such a subfactor $N\subset M$ exists, then
$[M_{k-1}:N]=p$ and $\mbox{depth}(N\subset M_{k-1})=2$
by Proposition~\ref{reducing depth}.
Therefore, $(N\subset M_{k-1})\cong (N\subset  N\rtimes \Bbb{Z}/p\Bbb{Z})$
(see \cite{NV2}, Corollary~4.19) and Theorem~\ref{galois correspondence}
implies the existence of a subgroup of $\Bbb{Z}/p\Bbb{Z}$
corresponding to the intermediate subfactor $N\subset M$.
But $\Bbb{Z}/p\Bbb{Z}$ does not have any non-trivial subgroups,
therefore $M=M_{k-1}$, i.e., $k=1$ and $n=2$.
\end{proof}
\medskip

Note that intermediate von Neumann subalgebras
of $M\subset M_1$ can be characterized in terms of projections in
$M'\cap M_2$ having certain properties \cite{B1}.
Namely, every projection $ q\in M'\cap M_2$ such that
\begin{enumerate}
\item[ (IS 1)]
$qe_2 =e_2$,
\item[ (IS 2)]
$E_{M_1}(q)$ is a scalar,
\item[ (IS 3)]
$\lambda^{-1}E_{M_1}(qe_1e_2)$ is a multiple of a projection,
\end{enumerate}
implements a conditional expectation from $M_1$ to an intermediate
subalgebra $Q=\{q\}'\cap M_1$. If $Q$ is a factor, then $[M_1:Q]
=\tau(q)^{-1}$
(\cite{B1}, Theorem 3.2). This result is true for all finite index
subfactors (regardless of depth).
The goal of next two propositions is to relate such projections and
coideal $*$-subalgebras in the case of finite depth inclusions
(which are intermediate subfactors of depth $2$ inclusions by
Corollary~\ref{subinclusions}).

\begin{proposition}
\label{introducing p_I}
If $I\subset B$ is a connected left coideal $*$-subalgebra, then
there is a constant $\lambda_I>0$ such that $p_I = \lambda_I^{-1}
H^{-1/2}x_IH^{-1/2}$ is a projection in $I$.
\end{proposition}
\begin{proof}
From the formula for $\Delta(x_I)$ (Proposition~\ref{x_I}) we get
$$
\eps_s(H^{-1}x_I) = m(S\otimes \id)\Delta(H^{-1}x_I)
= \sum_{\alpha rs}\, \frac{1}{\tau(f_{ss}^\alpha)} \,
f_{sr}^\alpha H^{-1} f_{rs}^\alpha \in Z(I)\cap B_s.
$$
Since $I$ is connected, we conclude that $\eps_s(H^{-1}x_I)
=\lambda_I 1$
for some constant $\lambda_I$.
Using this result and Proposition~\ref{x_I} one can check by a
direct computation that $\Delta(H^{-1}x_I)^2 = \lambda_I
\Delta(H^{-1}x_I)$, from where it follows that $\lambda_I^{-1}H^{-1}x_I$ is an
idempotent and $p_I = \lambda_I^{-1}H^{-1/2}x_IH^{-1/2}$ is a projection.
\end{proof}

\begin{proposition}
\label{IS for p_I}
Let $N\subset M$ be a depth $2$ inclusion, $B$  be the quantum groupoid
constructed on $M'\cap M_2$, and  $\theta :
M_1\rtimes B \to M_2$ be the isomorphism of II${}_1$ factors (see Preliminaries).
Then for any connected left coideal $*$-subalgebra $I\subset B$ the
projection $q_I =\theta(p_I)$ satisfies properties (IS 1)--(IS 3).
\end{proposition}
\begin{proof}
The relation $\eps_t(x_I H^{-1}) = \lambda_I 1$
follows easily from definitions, from where
$$
\eps_t(H^{1/2} p_I H^{-1/2}) = \lambda_I \eps_t(x_I H^{-1}) =1,
$$
hence $H^{1/2} p_I H^{-1/2}e_B = e_B$. Applying $\theta$ to
the last equality, and observing that $e_2 = \theta(H^{-1/2}e_BH^{-1/2})$,
we get $\theta(p_I)e_2 =e_2$ which is (IS 1).

To establish the second property, note that $E_{M_1}(\theta(p_I)) =
\theta(E_{M_1}(p_I)) = \theta(E_{B_t}(p_I))$, so it is enough to
show that $E_{B_t}(p_I)$ is a scalar. Using Proposition~\ref{expectation EI}
($x_{B_t} =H$) and Corollary~\ref{tau-invariance} we have :
\begin{eqnarray*}
E_{B_t}(p_I) &=& {p_I}\1 \tau(H{p_I}\2) =
\lambda_I^{-1} S(1\1) \tau(1\2x_I) = \lambda_I^{-1} 1.
\end{eqnarray*}

For the third property we must verify that $z_I = E_{M_1}(q_Ie_1e_2)$
is a multiple of a projection. Using the duality between $B^* =N'\cap M_1$
and $B=M'\cap M_2$ and the formula for $E_I$ from
Proposition~\ref{expectation EI} we have, for all $b\in B$ :
\begin{eqnarray*}
\la z_I,\, b \ra
&=&  \lambda^{-3} \tau(E_{M_1}(q_Ie_1e_2) e_2e_1Hb) \\
&=&  \tau(q_IHb) = \lambda_I^{-1} \tau(x_I G^{1/2}bG^{-1/2}),\\
&=& \lambda_I^{-1} \tau(x_I E_{I}(G^{1/2}bG^{-1/2})) \\
&=& \lambda_I^{-1} \tau(x_I (G^{1/2}bG^{-1/2})\1)
                   \tau(x_I (G^{1/2}bG^{-1/2})\2) \\
&=& \lambda_I \la z_I,\, b\1\ra \la z_I,\, b\2 \ra
    = \lambda_I \la z_I^2,\, b \ra,
\end{eqnarray*}
which implies that $z_I^2 = \lambda_I^{-1} z_I$.
Finally, in order to show that $z_I$ is selfadjoint,
we first compute $S(q_I)$, using properties of the antipode
and Proposition~\ref{x_I} :
\begin{eqnarray*}
S(q_I) &=& \lambda_I^{-1} S(S(H)^{1/2}H^{-1/2} x_I H^{-1/2} S(H)^{-1/2}) \\
       &=& \lambda_I^{-1} S(H)^{-1/2} H^{-1/2} x_I H^{-1/2} S(H)^{1/2}
             = S(H)^{-1} q_I S(H),
\end{eqnarray*}
from where, using the definition of $S$ we have
\begin{eqnarray*}
z_I^* &=& E_{M_1}(q_Ie_1e_2)^* = E_{M_1}(e_2e_1 S(H^{-1}q_I H))^* \\
      &=& E_{M_1}(e_2e_1q_I)^* = E_{M_1}(q_Ie_1e_2) = z_I.
\end{eqnarray*}
\end{proof}

\begin{remark}
Let $Q_I = \{q_I\}' \cap M_1$. Then $[M_1: Q_I] =\tau(q_I)
=\lambda_I^{-1}$ (cf.\ \cite{B1}).
\end{remark}

\begin{proposition}
\label{basic construction for I}
If $I$ is a left coideal $*$-subalgebra of $B$ and $\delta(I)$
is a left coideal subalgebra of $B^*$ constructed in
Proposition~\ref{dual lattices}, then the triple
$$
\delta(I) \subset B^* \subset B^*\rtimes I
$$
is a basic construction.
\end{proposition}
\begin{proof}
It follows from Proposition~\ref{IS for p_I} that $q_I$ implements
the conditional expectation from $N'\cap M_1=B^*$ to
$\{ q_I \}' \cap (N'\cap M_1) = \{ q_I \}' \cap B^* = \{ p_I \}' \cap B^*
= \tilde{I}' \cap B^* =\delta(I)$ (observe that $\tilde{I}$ is the
right coideal $*$-subalgebra of $B^*$ generated by $p_I$, cf.\
Propositions~\ref{tilde I} and \ref{introducing p_I}).
Since $B^*$ and $p_I$ generate $B^*\rtimes I$ we conclude that
$\delta(I) \subset B^* \subset B^*\rtimes I$ is a basic construction.
\end{proof}

\end{section}


\begin{section}
{Bimodule categories and principal graphs}

In this section we establish an equivalence between the tensor category
of ${N-N}$ bimodules of a finite index and finite depth subfactor $N\subset M$
and the co-representation category of a quantum groupoid $B$
canonically  associated with it as in Theorem~\ref{galois correspondence}.
The principal graph of $N\subset M$ can be described in terms of relative
Hopf modules over $B$. Alternatively, we show that it can also
be obtained as a certain Bratteli diagram.

Our methods follow those of \cite{KY},
where the special case of the invariants associated with the
group-subgroup subfactors was considered.

A left (resp., right) $B$-comodule $V$ (with the structure map
denoted by $v \mapsto v\I \otimes v\II,\, v\in V$)
is said to be {\em unitary}, if
\begin{eqnarray*}
(v_2\I)^*(v_1, v_2\II) &=& S(v_1\I)G (v_1\II,v_2) \\
(\mbox{resp., } (v_1\II)(v_1\I, v_2) &=& G^{-1}S((v_1\I)^*)(v_1,v_2\I)),
\end{eqnarray*}
where $v_1,v_2\in V,$ and $G$ is the canonical group-like element of
$B$. The notion of a unitary comodule in the Hopf
$*$-algebra case can be found, e.g., in (\cite{KS}, 1.3.2).

Given left coideal $*$-subalgebras $H$ and $K$ of $B$, we consider
a category $C_{H-K}$ of left relative $(B, H-K)$ Hopf bimodules
(cf. \cite{Takeuchi}), whose objects are Hilbert spaces which
are both $H-K$-bimodules and left unitary $B$-comodules such that the bimodule
action commutes with the coaction of $B$, i.e., for
any object $V$ of $C_{H-K}$ and $v\in V$ one has
$$
(h\lact v \ract k)\I \otimes (h\lact v \ract k)\II =
h\1 v\I k\1 \otimes (h\2 \lact v\II \ract k\2),
$$
where $v\mapsto v\I \otimes v\II$ denotes the coaction of $B$ on
$v$, $h\in H,\ k\in K$, and morphisms are  intertwining maps.

Similarly one can define a category of right relative
$(B, H-K)$ Hopf bimodules.

\begin{remark}
\label{on relative bimodules}
Note that any left $B$-comodule $V$ is automatically a $B_t-B_t$-bimodule via
$z_1\cdot v\cdot z_2 =  \eps(z_1v\I z_2)v\II,~v\in V, z_1,z_2\in B_t$.
For any object of $C_{H-K}$, this $B_t-B_t$-bimodule structure
is a restriction of the given $H-K$-bimodule structure; it is easily seen
by applying $(\eps\otimes id)$ to both sides of the relation of commutation
between the $H-K$-bimodule action and the coaction of $B$,
and taking $h,k\in B_t$.
Therefore, in the case when $H=B_t$ (resp.\ $K=B_t$)
we can speak about right (resp. \ left) relative Hopf modules,
which are a special case of weak Doi-Hopf modules \cite{Bohm}.

Let us also mention obvious relations
$\eps_s(v^{(1)})\otimes v^{(2)}= 1_{(1)}\otimes(v\ract 1_{(2)})$
and $\eps_t(v^{(1)})\lact v^{(2)}=v$.
\end{remark}

\begin{proposition}
\label{double cosets}
If $B$ is a  group Hopf $C^*$-algebra and $H,K$ are subgroups,
then there is a bijection between simple objects of $C_{H-K}$
and double cosets of $H\backslash B/K$.
\end{proposition}
\begin{proof}
If $V$ is an object of $C_{H-K}$, then every simple subcomodule
of $V$ is $1$-dimensional. Let $U =Cu$ ($u\mapsto g\otimes u,\, g\in B$)
be one of these comodules, then all other simple subcomodules of $V$
are of the form $h\lact U \ract k$, where $h\in H,\,k\in K$, and
$$
V= \oplus_{h,k}\, (h\lact U \ract k) = \span\{ HgK \}.
$$
Vice versa, $\span\{ HgK \}$ with natural $H-K$ bimodule and $B$-comodule
structures is a simple object of $C_{H-K}$.
\end{proof}
\begin{example}
\label{basic object}
If $H,V,K$ are left coideal $*$-subalgebras of $B,\ H\subset V,K\subset V$,
then $V$ is an object of $C_{H-K}$ with the structure maps given by
$ h\lact v \ract k = hvk$ and $\Delta$, where
$v\I \otimes v\II =\Delta(v),\ v\in V,\ h\in H,\ k\in K$. The scalar product is
defined by the restriction on $V$ of the Markov trace of $B$
(this $B$-comodule is unitary due to Corollary \ref{tau-invariance}).

Similarly, right coideal $*$-subalgebras of $B$ give examples of right relative
$(B, H-K)$ Hopf bimodules.
\end{example}

\medskip
Given an object $V$ of $C_{H-K}$, the conjugate Hilbert space
$\overline{V}$ is an object of $C_{K-H}$ with the bimodule action
$$
k\lact\overline{v}\ract h =\overline{h^*\lact v \ract k^*}\ \ \ (\forall h\in H,\, k\in K)
$$
(here $\overline{v}$ denotes the vector $v\in V$ considered as an
element of $\overline{V}$) and the coaction $\overline{v}\mapsto \overline{v}\I
\otimes \overline{v}\II=(v^{(1)})^*\otimes \overline{v}\II$. The relation of
commutation between the actions and the coaction and the unitarity of $\overline{V}$
are straightforward.

Define $V^*$, the dual object of $V$, to be $\overline{V}$ with
the above structures.
One can directly check that $V^{**}\cong V$ for any object  $V$.
Let us remark, that in Example \ref{basic object} the dual object can
be obtained by putting $\overline{v}=v^*$ for all $v\in V$.
\medskip

\begin{definition}
\label{reltenspr}
Let $L$ be another coideal $*$-subalgebra of $B$. For any objects
$V\in C_{H-L}$ and $W\in C_{L-K}$, we define an object $V\otimes_L W$
from $C_{H-K}$ as a tensor product of bimodules $V$ and $W$ \cite{JS}
equipped with a comodule structure
$$
(v\otimes_L w)^{(1)}\otimes (v\otimes_L w)^{(2)}:=v^{(1)}w^{(1)}\otimes
(v^{(2)}\otimes_L w^{(2)}).
$$
\end{definition}

Let us verify that we have indeed an object from $C_{H-K}$. First, the above
coproduct is clearly coassociative and compatible with counit
(see the properties of $\eps$ and Remark \ref{on relative bimodules}):
\begin{eqnarray*}
\eps(v^{(1)}w^{(1)})(v^{(2)}\otimes_L w^{(2)})
&=&  \eps(v^{(1)}1_{(1)})v^{(2)}\otimes_L \eps(1_{(2)}w^{(1)})w^{(2)} \\
&=&  (v\ract S(1_{(1)}))\otimes_L (1_{(2)}\lact w) = v\otimes_L w.
\end{eqnarray*}
Second, the commutation relation between the $H-K$-bimodule and
$B$-comodule structures can be proved as follows:
\begin{eqnarray*}
\lefteqn{ (h\lact (v\otimes_L w)\ract k)\I
\otimes (h\lact (v\otimes_L w)\ract k)\II = } \\
&=& ((h\lact v)\otimes_L (w\ract k))^{(1)}\otimes
    ((h\lact v)\otimes_L (w\ract k))^{(2)} \\
&=& (h\lact v)^{(1)}(w\ract k)^{(1)}\otimes (h\lact v)^{(2)}\otimes_L
    (w\ract k)^{(2)} \\
&=& h_{(1)}v^{(1)}w^{(1)}k_{(1)}\otimes (h_{(2)}\lact v^{(2)})\otimes_L
    (w^{(2)}\ract k_{(2)}) \\
&=& h_{(1)}(v\otimes_L w)^{(1)}k_{(1)}
    \otimes ( h_{(2)}\lact(v\otimes_L w)^{(2)}\ract k_{(2)} ).
\end{eqnarray*}
Finally, let us show that $V\otimes_L W$ is unitary.
To this end, recall the following expression of the
scalar product in this bimodule \cite{JS}:
\begin{eqnarray*}
\lefteqn{ (v_1\otimes_L w_1,v_2\otimes_L w_2)_{V\otimes_L W}
=(v_1\ract\la w_1,w_2\ra_L,v_2)_V= } \\
&=& (\la v_1,v_2\ra_L\lact w_1,w_2)_W = \tau(\la v_1,v_2\ra_L\la w_1,w_2\ra_L),
\end{eqnarray*}
where $v_1,v_2\in V, w_1,w_2\in W$ and the elements
$\la v_1,v_2\ra_L,\la w_1,w_2\ra_L\in L$ (~$L$-valued scalar products
on $V$ and $W$ respectively) are  defined in a unique way by the relations
$$
(v_1\ract l,v_2)_V=\tau(l\la v_1,v_2\ra_L),\qquad
(l\lact w_1,w_2)_W =\tau(l\la w_1,w_2\ra_L)\quad (\forall l\in L).
$$
Then the needed relation follows from
\begin{eqnarray*}
\lefteqn{ ((v_2\otimes_L w_2)^{(1)})^*
(v_1\otimes_L w_1,(v_2\otimes_L w_2)^{(2)})_{V\otimes_L W} =} \\
&=& (w_2^{(1)})^*(v_2^{(1)})^*(v_1\otimes_L w_1,v_2^{(2)}
\otimes_L w_2^{(2)})_{V\otimes_L W} \\
&=& (w_2^{(1)})^*(v_2^{(1)})^*(v_1\ract\la w_1,w_2^{(2)}\ra_L,v_2^{(2)})_V \\
&=& (w_2^{(1)})^*S(\la w_1,w_2^{(2)}\ra_L^{(1)})S(v_1^{(1)})G(v_1^{(2)}\ract
\la w_1,w_2^{(2)}\ra_L^{(2)},v_2)_V \\
&=& (w_2^{(1)})^*S(\la w_1,w_2^{(2)}\ra_L^{(1)})
\tau(\la w_1,w_2^{(2)}\ra_L^{(2)}\la v_1^{(2)},v_2\ra_L) S(v_1^{(1)})G \\
&=& (w_2^{(1)})^*(\la v_1^{(2)},v_2\ra_L^{(1)})^*G^{-1}
\tau(\la w_1,w_2^{(2)}\ra_L\la v_1^{(2)},v_2\ra_L^{(2)}) S(v_1^{(1)})G \\
&=& (w_2^{(1)})^*(\la v_1^{(2)},v_2\ra_L^{(1)})^*G^{-1}
 (\la v_1^{(2)},v_2\ra_L\lact w_1,w_2)_W^{(2)}  S(v_1^{(1)})G \\
&=& S(w_1^{(1)})(\la v_1^{(2)},v_2\ra_L\lact w_1^{(2)},w_2)_W S(v_1^{(1)})G \\
&=& S(w_1^{(1)})S(v_1^{(1)})G
(v_1^{(2)}\otimes_L w_1^{(2)},v_2\otimes_L w_2)_{V\otimes_L W} \\
&=& S((v_1\otimes_L w_1)^{(1)})
    G((v_1\otimes_L w_1)^{(2)},v_2\otimes_L w_2)_{V\otimes_L W},
\end{eqnarray*}
where we used the unitarity of $V$ and $W$ and Corollary~\ref{tau-invariance}.

\begin{lemma}[cf. \cite{KY}]
\label{properties of tensor product}
The operation of tensor product is
\newline (i)
associative, i.e., $V\otimes_L (W \otimes_P U)\cong(V\otimes_L W)\otimes_P U$,
\newline (ii)
compatible with duality, i.e., $(V\otimes_L W)^* \cong W^*\otimes_L V^*$,
\newline (iii)
distributive, i.e.,
$(V\oplus V')\otimes_L W = (V\otimes_L W) \oplus (V'\otimes_L W)$.
\end{lemma}
\begin{proof}
Easy excercise left to the reader.
\end{proof}

\medskip
The tensor product of morphisms $T\in \Hom(V, V')$ and $S\in
\Hom(W, W')$ is defined as usual :
$$
(T\otimes_L S)(v\otimes_L w) = T(v)\otimes_L S(w).
$$
From now on let us suppose that $B$ is biconnected and
acts outerly on the left on a II${}_1$ factor $N$
and that the extension of this action on $L^2(N)$ satisfies
the relations mentioned in the end of Preliminaries.

Given an object $V$ of $C_{H-K}$, we construct an
$N\rtimes H-N\rtimes K$-bimodule $\widehat{V}$ as follows.  We put
$$
\widehat{V} = \span \{ \Delta(1)\lact (\xi\otimes v)\mid
\xi\otimes v\in L^2(N)\otimes V\}
$$
and  denote $[\xi\otimes v] = \Delta(1)\lact(\xi\otimes v)$.
It is straightforward  to show that $\widehat{V}$ is
characterized by the property
$[\xi\otimes(z\lact v)]=[(\xi\ract z)\otimes v]
=[\xi(z\lact 1)\otimes v]\ \forall z\in B_t$, i.e., that
$\widehat{V} = L^2(N)\otimes_{B_t} V.$

Let us consider $\widehat{V}$ as a Hilbert space with the scalar product
$$
([\xi\otimes v],[\eta\otimes w])_{\widehat{V}}
=(\xi,\eta)_{\widetilde{L^2(N)}}(v,w)_V.
$$
and define the actions of $N,H,K$ on $\widehat{V}$ by
$$
a[\xi\otimes v] = [a\xi\otimes v], \qquad
[\xi\otimes v]a = [\xi(v\I\lact a)\otimes v\II],
$$
$$
h[\xi\otimes v] = [(h\1\lact \xi)\otimes (h\2\lact v)], \qquad
[\xi\otimes v]k = [\xi\otimes (v\ract k)],
$$
for all $a\in N,\, h\in H,\, k\in K$.  One can check that these
actions are well-defined and that
$$
ha[\xi\otimes v] = (h\1\lact a)h\2 [\xi\otimes {v}], \qquad
[\xi\otimes v]ka =  [\xi\otimes v](k\1\lact a)k\2,
$$
$$
(ha[\xi\otimes v])ka' = ha([\xi\otimes v]ka'),\qquad  a'\in N,
$$
i.e., that the above formulas define the structure of an
$(N\rtimes H)-(N\rtimes K)$ bimodule on $\widehat{V}$ in
the algebraic sense. Let us show that this bimodule
is unitary. We only need to check that
$$
([\xi\otimes v]a,[\eta\otimes w])_{\widehat{V}}=
([\xi\otimes v],[\eta\otimes w]a^*)_{\widehat{V}}
$$
(all other relations of unitarity are trivial). The following
computation uses the above definitions, the properties of the action
of $B$ on ${L^2(N)}$ and the unitarity of $V$:
\begin{eqnarray*}
([\xi\otimes v]a,[\eta\otimes v])_{\widehat{V}}
&=& ([\xi(v\I\lact a)\otimes v\II],[\eta\otimes v])_{\widehat{V}}\\
&=& (\xi(v\I\lact a),\eta)_{\widetilde{L^2(N)}}(v\II,w)_V\\
&=& (\xi(S(H)v\I\lact a)),\eta)_{L^2(N)}(v\II,w)_V\\
&=& (\xi,\eta(HS^{-1}((v\I)^*)\lact a^*))_{L^2(N)}(v\II,w)_V\\
&=& (\xi,\eta(GS^{-1}((v\I)^*)\lact a^*))_{\widetilde{L^2(N)}}(v\II,w)_V\\
&=& (\xi,\eta(w\I\lact a^*))_{\widetilde{L^2(N)}}(v,w\II)_V\\
&=& ([\xi\otimes v],[\eta\otimes w]a^*)_{\widehat{V}}.
\end{eqnarray*}
For any morphism $T\in \Hom(V,W)$, define a morphism
$\widehat{T} \in \Hom(\widehat{V},\widehat{W})$ by
$$
\widehat{T}([\xi\otimes v]) = [\xi\otimes T(v)].
$$
\begin{example}
\label{crossproduct as bimodule}
For $V\in C_{H-K}$ from Example \ref{basic object}, we have
$\widehat{V}=N\rtimes V$ as $N\rtimes H-N\rtimes K$-bimodules.
Indeed, the algebraic operations and the scalar products are the same
(see the definitions above and Remark~\ref{Markov cond.exp}).
\end{example}

\begin{theorem}
\label{properties of the functor}
The above assignments $V\mapsto \widehat{V}$ and $T\mapsto \widehat{T}$
define a functor from  $C_{H-K}$ to the category of
$N\rtimes H - N\rtimes K$ bimodules. This functor preserves direct
sums and is compatible with operations of taking tensor products and
adjoints in the sense that if $W$ is an object of $C_{K-L}$, then
$$
\widehat {V\otimes_K W} \cong \widehat{V} \otimes_{N\smrtimes K} \widehat{W},
\quad \mbox{ and } \quad \widehat{V^*} \cong (\widehat{V})^*.
$$
\end{theorem}
\begin{proof}
1. Directly from definitions we have $\widehat{V\oplus W} =
\widehat{V} \oplus \widehat{W}$ for any objects
$V$ and $W$ of $C_{H-L}$.
\newline 2. In order to show that $\widehat{V\otimes_L W} \cong
\widehat{V}\otimes_{N\smrtimes L} \widehat{W}$, let us define a map
$$
\alpha : \widehat{V\otimes_L W } \to
\widehat{V} \otimes_{N\smrtimes L} \widehat{W} :
[\xi \otimes (v\otimes_L w)] \mapsto  [\xi \otimes v]
\otimes_{N\smrtimes L} [1\otimes w],
$$
for all $\xi\in L^2(M),\,v\in V,\,w\in W$.
\newline
(a) $\alpha$ is well-defined, since, for all
$\xi\in L^2(N),  z\in B_t, v\in V, w\in W$
\begin{eqnarray*}
\alpha([(\xi\cdot z) \otimes (v\otimes_L w)])
&=& [\xi \otimes \eps(zv\I)v\II] \otimes_{N \smrtimes L} [1\otimes w] \\
&=& [\xi \otimes \eps(zv\I1\1)(v\II\cdot 1\2)] \otimes_{N \smrtimes L}
    [1\otimes w] \\
&=& [\xi \otimes \eps(zv\I1\1)v\II] \otimes_{N \smrtimes L}
    [(1\2 \cdot 1)\otimes w] \\
&=& [\xi \otimes \eps(zv\I1\1)v\II] \otimes_{N \smrtimes L}
    [1\otimes \eps(1\2w\I)w\II] \\
&=& [\xi \otimes \eps(zv\I w\I)v\II] \otimes_{N \smrtimes L} [1\otimes w\II] \\
&=& \alpha([\xi \otimes \eps(zv\I w\I)(v\II \otimes_L w\II)]) \\
&=& \alpha([\xi \otimes z\cdot (v\otimes_L w)]).
\end{eqnarray*}
(b) $\alpha$ preserves the bimodule structure.
Indeed, for all $h\in H, k\in K, a,a' \in N$ we have :
\begin{eqnarray*}
\alpha(ah[\xi \otimes (v\otimes_L w)])
&=& \alpha([a(h\1\lact \xi) \otimes ((h\2\lact v) \otimes_L w)]) \\
&=& [a(h\1\lact \xi) \otimes (h\2\lact v)] \otimes_{N \smrtimes L}
    [1\otimes w] \\
&=& ah [\xi \otimes v] \otimes_{N \smrtimes L} [1\otimes w] \\
&=& ah \alpha([\xi \otimes(v\otimes_L w)]), \\
\alpha([\xi \otimes (v\otimes_L w)]a'k)
&=& \alpha([\xi(v\I w\I\cdot a') \otimes (v\II \otimes_L (w\II\ract k))]) \\
&=& [\xi(v\I w\I\cdot a') \otimes v\II] \otimes_{N \smrtimes L}
    [1\otimes (w\II\ract k)] \\
&=& [\xi \otimes v] \otimes_{N \smrtimes L}
    [(w\I \cdot a') \otimes (w\II\ract k)] \\
&=& [\xi \otimes v] \otimes_{N \smrtimes L} ([1\otimes w]a'k) \\
&=& \alpha([\xi \otimes(v\otimes_L w)])a'k.
\end{eqnarray*}
(c) Observe that the map
$$
\beta : \widehat{V} \otimes_{N\smrtimes L} \widehat{W}
\to \widehat{V\otimes_L W } :
[\xi \otimes v] \otimes_{N\smrtimes L} [\eta\otimes w]\mapsto [\xi
(v\I\cdot\eta)\otimes(v\II\otimes_L w)],
$$
where $\xi,\eta\in L^2(N),v\in V,w\in W$, is the inverse of $\alpha$.
\newline  (d) $\alpha$ is an isometry of Hilbert spaces.
Indeed, by the above definitions,
$$
\vert\vert[\xi\otimes(v\otimes_L w)]\vert\vert^2_{\widehat{V\otimes_L W}}
=(S(H)\lact \xi,\xi)_{L^2(N)}(v\ract\la w,w\ra_L,v)_V,
$$
where the $L$-valued scalar product $\la w,w\ra_L$ on $W$ is defined in
a unique way by the relation $(l\lact w,w)_W=\tau(l\la w,w\ra_L)\
(\forall l\in L)$.

On the other hand, the definition of the tensor product of bimodules (see,
for example, \cite{JS}) gives:
$$
\vert\vert[\xi \otimes v]\otimes_{N\smrtimes L} [1\otimes w]\vert\vert^2_
{\widehat{V}\otimes_{N\smrtimes L} \widehat{W}}=
([\xi\otimes v]\ract\la [1\otimes w],[1\otimes w]\ra_{N\smrtimes L},
[\xi\otimes v])_{\hat V},
$$
where the element $\la [1\otimes w],[1\otimes w]\ra_{N\smrtimes L}$ is
defined in a unique way by
$$
([n\otimes l]\lact [1\otimes w],[1\otimes w])_{\hat W}
=\tau_{N\smrtimes L}([n\otimes l]
\la [1\otimes w],[1\otimes w]\ra_{N\smrtimes L}).
$$
Since the left-hand side of this equality can be rewritten as:
\begin{eqnarray*}
\lefteqn{ ([n(l_{(1)}\lact 1)\otimes
(l_{(2)}\lact w)],[1\otimes w])_{\hat W} =} \\
&=& ([n\otimes(\eps_t(l_{(1)})l_{(2)})\lact w)],[1\otimes w])_{\hat W}
\tau_N(S(H)\lact n)(l\lact w,w)_W,
\end{eqnarray*}
we can see that $\la [1\otimes w],[1\otimes w]\ra_{N\smrtimes L}=
[1\otimes\la w,w\ra_L]$. Together with the formula for the scalar product
on $\hat V$  this gives
$$
\vert\vert[\xi\otimes(v\otimes_L w)]\vert\vert^2_{\widehat{V\otimes_L W}}=
\vert\vert[\xi \otimes v]\otimes_{N\smrtimes L} [1\otimes w]\vert\vert^2_
{\widehat{V}\otimes_{N\smrtimes L} \widehat{W}}.
$$
3. In order to show that $\widehat{V^*} \cong (\widehat{V})^*$,
let us define a map
$$
\gamma:\widehat{V^*} \to (\widehat{V})^*:[\xi\otimes\overline v]
\mapsto\overline{[(v^{(1)}\lact J\xi)\otimes v^{(2)}]}.
$$
(a) $\gamma$ is well-defined, since, for all $\xi\in L^2(N),
z\in B_t, v\in V$
\begin{eqnarray*}
\lefteqn{ \gamma([(S(z)\lact\xi)\otimes\overline v])=
\overline{[(v^{(1)}\lact J(S(z)\lact\xi))\otimes v^{(2)}]} =} \\
&=& \overline{[(v^{(1)}z^*\lact \xi)\otimes v^{(2)}]}
    =\gamma([\xi\otimes(z\lact \overline v)]).
\end{eqnarray*}
(b) It is straightforward to show that $\gamma$
preserves the bimodule structure and that the map
$$
\overline{[\xi\otimes v]}\mapsto
[((v^{(1)})^*\lact J\xi)\otimes \overline v^{(2)}]
$$
from $(\widehat{V})^*$ to $\widehat{V^*}$ is the inverse of $\gamma$.
\newline(c) $\gamma$ is an isometry of Hilbert spaces.
Indeed, using the above definitions, the relation
$(J\xi,J\eta)_{\widetilde{L^2(M_1)}}=(G\lact \eta,\xi)_{\widetilde{L^2(M_1)}}$
and the unitarity of $V^*$, we have:
\begin{eqnarray*}
\vert\vert \gamma([\xi\otimes\overline v])\vert\vert_{(\widehat{V})^*}^2
&=&  \vert\vert \overline{[(v^{(1)}\lact J\xi)\otimes
v^{(2)}]}\vert\vert_{(\widehat{V})^*}^2 \\
&=& \vert\vert v^{(1)}\lact J\xi\vert\vert_{\widetilde{L^2(N)}}^2
\vert\vert\overline v^{(2)}\vert\vert_{V^*}^2 \\
&=& \vert\vert G^{1/2}S^{-1}(\overline v^{(1)})\lact
\xi\vert\vert_{\widetilde{L^2(N)}}^2
\vert\vert\overline v^{(2)}\vert\vert_{V^*}^2 \\
&=& (\xi,(\overline v^{(1)})^{(2)}
S^{-1}((\overline v^{(1)})^{(1)})\lact \xi)_{\widetilde{L^2(N)}}
(\overline v,\overline v^{(2)})_{V^*}\\
&=&
(\xi,S(\eps_t(\overline v^{(1)}))\lact \xi)_{\widetilde{L^2(N)}}
(\overline v,\overline v^{(2)})_{V^*} \\
&=& \vert\vert\xi\vert\vert_{\widetilde{L^2(N)}}
(\overline v,\eps_t(\overline v^{(1)})\lact\overline v^{(2)})_{V^*}
= \vert\vert [\xi\otimes\overline v]\vert\vert_{\widehat{V^*}}^2.
\end{eqnarray*}
\end{proof}

According to Remark~\ref{on relative bimodules},  $C_{B_t-B_t}$
is nothing but the category of $B$-comodules, $\Corep(B)$. The next
theorem shows that this category is equivalent to $Bimod_{N-N}(N\subset M)$,
the category of $N-N$ bimodules of a subfactor $N\subset M$.
Recall that the latter is the tensor category generated by
simple subobjects of ${}_N (M_n)_N,\, n\geq 1$.

\begin{theorem}
\label{N-N equivalence}
Let $N\subset M$ be a finite depth subfactor ($[M:N]<\infty$), $k$
be a number such that$N\subset M_k$ has depth $\leq 2$,
and let $B$ be a canonical quantum groupoid
such that $(N\subset M_k)\cong (N\subset N\rtimes B)$.
Then  $Bimod_{N-N}(N\subset M)$ and  $\Rep(B^*)$
are equivalent as tensor categories.
\end{theorem}
\begin{proof}
First, we observe that
$$
Bimod_{N-N}(N\subset M) = Bimod_{N-N}(N\subset M_l)
$$
for any $l\geq 0$. Indeed, since both categories are semisimple, it is
enough to check that they have the same set of simple objects.
Clearly, all objects
of $Bimod_{N-N}(N\subset M_l)$ are objects of
$Bimod_{N-N}(N\subset M)$. Conversely, since irreducible $N-N$ subbimodules
of ${}_N L^2(M_i)_N$ are contained in the decomposition of
${}_N L^2(M_{i+1})_N$ for all $i\geq 0$, we see that objects of
$Bimod_{N-N}(N\subset M)$ belong to $Bimod_{N-N}(N\subset M_l)$.

Hence, by Proposition~\ref{reducing depth},
the problem can be reduced to the case when $N\subset M$ has
depth $2$ ($M=N\rtimes B$), i.e., it will suffice to prove that
$Bimod_{N-N}(N\subset N\rtimes B)$ is equivalent to $\Corep(B)$.
The previous theorem gives a functor from
$\Corep(B) = \Rep(B^*)$ to $Bimod_{N-N}(N\subset N\rtimes B)$.
To prove that this functor is, in fact, an equivalence,
we  need to check that it yields a bijection
between classes of simple objects of these categories.

Observe that $B$ itself is an object of $\Corep(B)$ via
$\Delta : B\to B\otimes B$ and $\widehat{B} = {}_N L^2(M)_N$.
Since the inclusion $N\subset M$ has depth $2$, the simple
objects of $Bimod_{N-N}(N\subset M)$
are precisely irreducible subbimodules of ${}_N L^2(M)_N$.
We have $\widehat{B} = {}_N L^2(M)_N = \oplus_i\,   {}_N
(p_iL^2(M))_N$, where $\{p_i\}$ is a family of mutually orthogonal
minimal projections in $N'\cap M_1$ so that every
bimodule $p_iL^2(M)$ is irreducible. On the other hand, $B$ is
cosemisimple, hence $B=\oplus_i\, V_i$, where each $V_i$
is an irreducible subcomodule. Note that $N'\cap M_1= B^* =\sum p_iB^*$
and every $p_iB$ is a simple submodule of $B$ ($=$ simple subcomodule
of $B^*$). Thus, we see that there is a bijection between simple
objects of $\Corep(B)$ and $Bimod_{N-N}(N\subset M)$, so that
the categories are equivalent.
\end{proof}
\medskip

The principal and dual principal graphs of a subfactor $N\subset M$
are defined as follows \cite{JS}, \cite{GHJ}, \cite{KY}.
Let $X= {}_N L^2(M)_M$ and consider the following sequence of
$N-N$ and $N-M$ bimodules :
$$
{}_N L^2(N)_N,~X,~X\otimes_M X^*,~X\otimes_M X^*\otimes_N X,\dots
$$
obtained by right tensoring with $X^*$ and $X$. The vertex set
of the principal graph is indexed by the classes of simple bimodules
appearing as summands in the the above  sequence. We connect
vertices corresponding to bimodules ${}_N Y_N$ and ${}_N Z_M$ by $l$
edges if ${}_N Y_N$ is contained in the decomposition of ${}_N Z_N$,
the restriction of ${}_N Z_M$, with multiplicity $l$.

The dual principal graph can be constructed in a similar way
from the following $M-M$ and $M-N$ bimodules :
$$
{}_M L^2(M)_M,~X^*,~X^*\otimes_N X,~X^*\otimes_N X\otimes_M X^*,\dots.
$$

We apply Theorem~\ref{properties of the functor} to express the principal
and dual graphs of a finite depth subfactor $N\subset M$ in terms of the
quantum groupoid associated with it.

Let $B$ and $K$ be a quantum groupoid and its left coideal
$*$-subalgebra such that $B$ acts on $N$ and $(N\subset M) \cong
(N\subset N\rtimes K)$. Then $\widehat{K} = {}_N L^2(M)_M$, where
we view $K$ as a relative $(B, K)$ Hopf module as in
Example~\ref{basic object}

By Theorem~\ref{properties of the functor} we can identify
irreducible $N-N$ (resp.\ $N-M$) bimodules with  simple $B$-comodules
(resp.\  relative right $(B, K)$ Hopf modules). Consider a bipartite
graph with vertex set given by the union of (classes of) simple
$B$-comodules and simple relative right $(B, K)$ Hopf modules
and the number of edges between the vertices $U$ and $V$ representing
$B$-comodule and relative right $(B, K)$ Hopf module respectively
being  equal to the multiplicity of $U$ in the decomposition of $V$
(when the latter is viewed as a $B$-comodule) :
\begin{eqnarray*}
& \mbox{simple}~B\mbox{-comodules}& \\
& \cdots \vert \cdots \vert \cdots & \\
& \mbox{simple relative right}~(B,K)~\mbox{Hopf modules} &
\end{eqnarray*}
The principal graph of $N\subset M$ is the connected part of the above graph
containing the trivial $B$-comodule.

Similarly, the dual principal graph can be obtained from
the following diagram
\begin{eqnarray*}
& \mbox{simple relative}~(B, K-K)~\mbox{Hopf bimodules} & \\
& \cdots \vert \cdots \vert \cdots & \\
& \mbox{simple relative left}~(B,K)~\mbox{Hopf modules} &
\end{eqnarray*}
as the connected component containing the relative Hopf $(B, K-K)$
bimodule $K$ (it corresponds to the bimodule ${}_M L^2(M)_M$.

Using the antiisomorphism $K \mapsto \delta(K)$ between
the lattices of left coideal subalgebras of $B$ and $B^*$
from Proposition~\ref{dual lattices}, it is possible to
express the principal graph of $N\subset N\rtimes K$ as a certain
Bratteli diagram.

\begin{proposition}
\label{Bratteli diagram = principal graph}
If $K$ is a coideal $*$-subalgebra of $B$ then the principal
graph of the subfactor $N\subset N\rtimes K$ is given by
the connected component of the Bratteli
diagram of the inclusion $\delta(K) \subset B^*$
containing the trivial representation of $B^*$.
\end{proposition}
\begin{proof}
First, let us show that there is a bijective correspondence
between right relative $(B,K)$ Hopf modules and $(B^*\rtimes K)$-modules.
Indeed, every right $(B,K)$ Hopf module $V$ carries a right action of $K$.
If we define a right action of $B^*$ by
$$
v\ract x = \la v\I,\, x\ra v\II,\quad v\in V, x\in B^*,
$$
then we have
\begin{eqnarray*}
(v \ract k)\ract x
&=& \la v\I k\1,\,x \ra (v\II \ract k\2) \\
&=& \la v\I,\, (k\1\lact x) \ra (v\II \ract k\2) \\
&=& (v \ract (k\1\lact x)) \ract k\2,
\end{eqnarray*}
for all $x\in B^*$ and $k\in K$ which shows that $kx$
and $(k\1\lact x)k\2$ act on $V$ exactly in the same way,
therefore $V$ is a right $(B^*\rtimes K)$-module.

Conversely, given an action of $(B^*\rtimes K)$ on $V$,
we automatically have a $B$-comodule structure such that
\begin{eqnarray*}
\la v\I k\1,x \ra (v\II \lact k\2)
&=& \la v\I,\, x\1\ra \la k\1,\, x\2 \ra (v\II \ract k\2) \\
&=& (v\ract (k\1\lact x)) \ract k\2 = (v\ract k)\ract x \\
&=& \la x,\, (v \ract k)\I \ra (v \ract k)\II,
\end{eqnarray*}
which shows that $v\I k\1 \otimes (v\II \ract k\2) =
(v \ract k)\I \otimes (v \ract k)\II$, i.e., that $V$
is a right relative $(B,K)$-module.

Thus, we see that the principal graph is given by
the connected component the Bratteli diagram
of the inclusion $B^* \subset B^*\rtimes K$
containing the trivial representation of $B^*$.
Recall that $B^*\rtimes K$
is the basic construction for the inclusion $\delta(K)\subset B^*$,
therefore the Bratteli diagrams of the above two inclusions are the same.
\end{proof}

\begin{corollary}
\label{principal graph for depth 2}
If $N\subset N\rtimes B$ is a depth 2 inclusion
corresponding to the quantum groupoid $B$, then
its principal graph is given by the Bratteli diagram
of the inclusion $B^*_t \subset B^*$.
\end{corollary}
\begin{proof}
In this case $K=B$ and inclusion $B^*_t \subset B^*$
is connected, so that $\delta(K) = B^*_t$ (note that
$B^*$ is biconnected).
\end{proof}

Let us mention two properties of the set $X_n$
of finite index values of subfactors with depth $\leq n$.

\begin{remark}
(a)
We can use Corollary~\ref{principal graph for depth 2}
to give a short proof of the fact that for any given
$n$ the set $X_n$ is a discrete subset of 
$\{4\cos^2 \frac{\pi}{n} \mid n\geq 3\}\cup [4,+\infty)$ \cite{J}. 

It follows from Corollary~\ref{subinclusions} that the index of any
depth  $\leq n$ subfactor is the $n$-th root of the index of
a depth $\leq 2$ subfactor, therefore we have
$X_n =\{ \root n \of x \mid x\in X_2\}$,
therefore it suffices to prove that $X_2$ is discrete.

Let $B$ be a biconnected quantum groupoid and
$\Lambda$ be the inclusion matrix of $B_t \subset B$.
We will show that all the entries of $\Lambda\Lambda^t$
are strictly positive. Indeed, let $\pi_1,\dots \pi_N$
(resp.\ $\rho_1,\dots \rho_M$) be all the classes of
irreducible representations of $B_t$ (resp.\ $B$), and
assume that $\rho_1$ is the trivial representation
of $B$ on $B_t$ (i.e., $\rho_1(b)z =\eps_t(bz)$ for all $b\in B,\, z\in B_t$
(\cite{BNSz}, 2.4, \cite{NV1}, 2.2).
Then $\Lambda_{ij}$, the $ij$-th entry of $\Lambda$,
is equal to the multiplicity of $\pi_i$ in $\rho_j\vert_{B_t}$.

Since $\rho_1\vert_{B_t}$ is faithful, we have $\Lambda_{i1}>0$
for all $i=1\,\dots M$, therefore
$$
(\Lambda\Lambda^t)_{ik} =\Sigma_j\,\Lambda_{ij}\Lambda_{kj}
\geq \Lambda_{i1}\Lambda_{k1} >0.
$$

Thus, it follows from Corollary~\ref{principal graph for depth 2}
that every element of $X_2$ is the norm of a matrix with strictly
positive entries. But for any given $m$ the number of
such matrices with norm $\leq m$ is clearly finite,  hence
$X_2 \cap (0,\,m]$ is finite for every $m$, i.e., $X_2$ is discrete.

(b) $X_n$ is also a multiplicative subsemigroup of $\Bbb{R}^{+}$.
Indeed, if $N\subset M$ and $P\subset Q$ are two subfactors
of depth $\leq n$ then $(N\otimes P) \subset (M\otimes Q)$
has depth $\leq n$ and $[(M\otimes Q):(N\otimes P)]= [M:N][Q:P]$.
\end{remark}

\end{section}


\begin{section}
*{Appendix : The structure of a quantum groupoid associated with a
finite depth subfactor}

We will write down explicit formulas that define
a quantum groupoid canonically associated with a finite
depth subfactor $N \subset M$ ($[M:N] =\lambda ^{-1}$).
It follows from Proposition~\ref{reducing depth} that
the subfactor $N \subset M_k$ is of depth $2$ for $k$
large enough. According to \cite{PP2}, the Jones tower
for the latter inclusion is
$$
N \subset M_k \subset M_{2k+1} \subset M_{3k+2} \subset \cdots
$$
Therefore, there is a non-degenerate duality between algebras
$A =N'\cap M_{2k+1}$ and $B=M_{k}'\cap M_{3k+2}$ making them
quantum groupoids dual to each other \cite{NV2}.
The corresponding bilinear form (cf.\ Preliminaries) is
given by
$$
 \la a,\, b\ra = \lambda^{-2(k+1)}\tau(af_2f_1Hb),\qquad a\in A, b\in B,
$$
where $H=\Index \tau\vert_{M_k'\cap M_{2k+1}}$ and
\begin{eqnarray*}
f_1 &=& \lambda^{k(k+1)/2}(e_{k+1}e_k\dots e_1)(e_{k+2}\dots e_2)
        \dots (e_{2k+1}\dots e_{k+1}), \\
f_2 &=& \lambda^{k(k+1)/2}(e_{2k+2}e_{2k+1}\dots e_{k+2})
         (e_{2k+3}\dots e_{k+3}) \dots (e_{3k+2}\dots e_{2k+2}),
\end{eqnarray*}
are the Jones projections of the $k$-step basic construction \cite{PP2}
such that $f_1$ (resp.\ $f_2$) implements the conditional expectation from
$M_{k}$ (resp.\ $M_{2k+1}$) to $N$ (resp.\ $M_k$).

The target and source counital subalgebras of $B$ are
$B_t = M_{k}'\cap M_{2k+1}$ and $B_s = M_{2k+1}'\cap M_{3k+2}$.
Note that $B$ is generated by $B_s$, $B_t$, and $e_{2k+2}$
as an algebra. Indeed, $M_{k}'\cap M_{3k+2} =
\la M_{k}'\cap M_{2k+1},\, e_{2k+2},\dots e_{3k+2} \ra$
because of the finite depth condition.

The antipode of $B$ is given by
$$
S(b) = j(HbH^{-1}), \quad b\in B,
$$
where $j(b) =J_{2k+2} b^* J_{2k+2}$ is a canonical $*$-anti-isomorphism
of $B = M_{k}'\cap M_{3k+2}$ (here $J_{2k+2}$ is the modular involution
on $L^2(M_{2k+1})$).

It is convenient to describe the comultiplication in terms of
separability elements. By definition, a {\em separability element}
(\cite{P}, 10.2) of a finite-dimensional $C^*$-algebra $D$ is a projection
$P_D\in D^{op}\otimes D$ uniquely determined by the properties
$$
(x_1\otimes 1)P_D(x_2\otimes 1) =(x_2\otimes 1)P_D(x_1\otimes 1),
\quad\mbox{ and }\quad m(P_D)=1,
$$
for all $x_1,x_2\in D$, where $m$ denotes the multiplication of $D$.

Let $I=M_{k}'\cap M_{2k+2}$, then we have
\begin{eqnarray*}
\Delta(yz) &=& (z\otimes y)\cdot(S\otimes \id)P_{B_t},\quad
               y\in B_s,\, z\in B_t, \\
\Delta(e_{2k+2}) &=& (S\otimes \id)P_I.
\end{eqnarray*}
Indeed, the first formula holds true in every weak Hopf algebra
(see, e.g., \cite{NV1}). To establish the second formula for all
$a,c\in A$ we compute (using the notation $P_I = P_I\I \otimes P_I\II$) :
\begin{eqnarray*}
\lefteqn{ \la a,\,H^{-1}S(P_I\I) \ra \la c,\, P_I\II \ra =}\\
&=& \lambda^{-4(k+1)}\,\tau(af_2f_1S(P_I\I H))\, \tau(cf_2f_1P_I\II) \\
&=& \lambda^{-4(k+1)}\,\tau(HP_I\I f_1f_2a)\, \tau(cf_2f_1 P_I\II) \\
&=& \lambda^{-4(k+1)+1}\,\tau(\lambda^{-1}HP_I\I E_{M_k'\cap M_{2k+2}}(f_1f_2a))
    \,\tau(cf_2f_1P_I\II) \\
&=& \lambda^{-4(k+1)+1}\, \tau(cf_2f_1 E_{M_k'\cap M_{2k+2}}(f_1f_2a)) \\
&=& \lambda^{-3(k+1)+1}\, \tau(cf_2f_1E_{M_{2k+2}}(f_1 E_{M_k'}(f_1f_2a))) \\
&=& \lambda^{-2(k+1)}\, \tau(cf_2f_1e_{2k+2}a) =\la ac,\, e_{2k+2}\ra ,
\end{eqnarray*}
where we used that $\Index\tau\vert_{M_k'\cap M_{2k+2}}=\lambda^{-1}H$
and that $E_{M_{2k+2}}(f_2) = \lambda^k e_{2k+2}$. Thus,
$\Delta(H^{-1}e_{2k+2}) = (H^{-1}\otimes 1)\cdot(S\otimes \id)P_I$
and, therefore, $\Delta(e_{2k+2}) = (S\otimes \id)P_I$.

Finally, the counit is given by
$$
\eps(b) =\lambda^{-(k+1)}\tau(f_2Hb), \quad b\in B.
$$

Note that $I = M_{k}'\cap M_{2k+2}$ is a left coideal $*$-subalgebra
in $B$ and that
$$
(N\subset M) \cong (M_{2k+1}\subset M_{2k+2})\cong (M_k\subset M_k\rtimes I).
$$

An example of a quantum groupoid of dimension $13$, associated to the
subfactor with index $4\cos^2\frac{\pi}{5}$,
was considered in (\cite{NV2}, 7.3). One can also describe
the quantum groupoids corresponding to the whole sequence
of subfactors with principal graphs $A_n$ \cite{J}. 

\end{section}



\begin{thebibliography}{99}

\bibitem{B1}
D.~Bisch.
\newblock A note on intermediate subfactors.
\newblock {\em Pacific J. Math.}, {\bf 163}, no.2 (1994), 201--216.

\bibitem{Bohm}
G.~Bohm.
\newblock Doi-Hopf modules over weak Hopf algebras.
\newblock {\em math.QA}/{\bf 9905027}(1999).

\bibitem{BNSz}
G.~B\"ohm, F.~Nill, and K.~Szlach\'anyi.
\newblock  Weak Hopf algebras I: Integral theory and $C^*$-structure.
\newblock {\em math.QA}/{\bf 9805116} (1998).

\bibitem{E}
M.~Enock.
\newblock Sous-facteurs interm\'ediaires et groupes quantiques mesur\'es.
\newblock {\em J. Operator Theory}, {\bf 42} (1999), no. 2, 305--330

\bibitem{GHJ}
F.M.~Goodman, P.~de la Harpe, and V.F.R.~Jones.
\newblock Coxeter graphs and towers of algebras. MSRI Publ. 14,
\newblock Springer-Verlag, (1989).

\bibitem{ILP}
M.~Izumi, R.~Longo, and S.~Popa.
\newblock A Galois correspondence for compact groups
of automorphisms of von Neumann algebras with a generalization to Kac algebras,
\newblock {\em J. Funct. Anal.}, {\bf 155} (1998), no. 1, 25--63.

\bibitem{KS}
A.~Klimyk and K.~Schm\"udgen.
\newblock Quantum Groups and Their Representations,
\newblock {\em Springer Verlag}, (1997).

\bibitem{KY}
H.~Kosaki and  S.~Yamagami.
\newblock Irreducible bimodules associated with crossed product algebras,
\newblock {\em Intern. J. Math.}, {\bf 3}(1992), 661-676.

\bibitem{J}
V.~Jones.
\newblock Index for subfactors.
\newblock {\em Invent. math.}, {\bf 72}(1983), 1-25.

\bibitem{JS}
V.~Jones and V.S.~Sunder.
\newblock Introduction to subfactors,
\newblock Cambridge University Press, (1997).

\bibitem{N2}
D.~Nikshych.
\newblock Duality for actions of weak Kac algebras and
crossed product inclusions of II${}_1~$factors.
\newblock {\em math.QA}/{\bf 9810049} (1998).

\bibitem{NV1}
D.~Nikshych and L.~Vainerman.
\newblock Algebraic versions of a finite dimensional quantum groupoid.
\newblock To appear in  {\em Lecture Notes in Pure and Appl. Math.},
\newblock {\em math.QA}/{\bf 9808054} (1998).

\bibitem{NV2}
D.~Nikshych and L.~Vainerman.
\newblock A characterization of depth 2 subfactors of II${}_1$ factors.
\newblock To appear in  {\em J. Func. Analysis},
\newblock {\em math.QA}/{\bf 9810028} (1998).

\bibitem{NSzW}
F.~Nill, K. Szlach\'anyi, and H.-W.~Wiesbrock.
\newblock Weak Hopf algebras and reducible Jones inclusions of
depth 2. I: From crossed products to Jones towers.
\newblock  {\em math.QA}/{\bf 9806130} (1998).

\bibitem{P}
R.~Pierce.
\newblock Associative algebras.
\newblock Graduate Texts in Mathematics, 88.
{\em Springer Verlag}, (1982).

\bibitem{PP2}
M.~Pimsner and S.~Popa.
\newblock Iterating the basic construction.
\newblock {\em Trans. Amer. Math. Soc.} {\bf 310} (1988), no. 1,
127--133.



\bibitem{Takeuchi}
M.~Takeuchi,
\newblock Relative Hopf modules -- equivalencies and freeness criteria.
\newblock {\em J. Algebra}, {\bf 60} (1979), 452-471.


\bibitem{Y}
T.~Yamanouchi.
\newblock Duality for generalized Kac algebras and a characterization
of finite groupoid algebras.
\newblock {\em J. Algebra}, {\bf 163} (1994), 9--50.

\bibitem{W1}
Y.~Watatani.
\newblock Index for $C^*$-subalgebras,
\newblock {\em Memoirs of the AMS} {\bf 424} (1990).

\bibitem{W}
Y.~Watatani.
\newblock Lattices of intermediate subfactors.
\newblock  {\em J. Funct. Anal.}, {\bf 140} (1996), no. 2, 312--334.

\end{thebibliography}
\end{document}